\documentclass[a4paper,11pt]{amsart}
\usepackage{amssymb}
\usepackage{amsmath}
\usepackage{mathtools}
\usepackage[T1]{fontenc}
\usepackage{amsthm}
\usepackage{amsfonts}
\usepackage{enumerate}
\usepackage{mathrsfs}
\usepackage{hyperref}
\usepackage{pdfpages}

\usepackage{pgfplots}
\pgfplotsset{compat=1.15}
\usetikzlibrary{arrows}

\usepackage{caption}
\usepackage[labelformat=simple]{subcaption}

\usepackage{xcolor}

\newlength\mylen

\hypersetup{
	colorlinks   = true, 
	urlcolor     = blue, 
	linkcolor    = blue, 
	citecolor   = red 
}

\usepackage[a4paper]{geometry}
\newtheorem{thm}{Theorem}[section]
\newtheorem{prop}[thm]{Proposition}
\newtheorem{lem}[thm]{Lemma}
\newtheorem{coro}[thm]{Corollary} 

\newtheorem*{question*}{Question}
\newtheorem*{thm*}{Theorem}
\theoremstyle{definition}

\theoremstyle{remark} 
\newtheorem{rem}[thm]{Remark}
\newtheorem*{rem*}{Remark}
\theoremstyle{definition}

\newtheorem*{nota*}{Notation}

\numberwithin{equation}{section}

\newcommand{\OK}{\mathcal O_{\mathbb K}}
\newcommand{\bP}{\boldsymbol P}
\newcommand{\X}{{\boldsymbol X}}
\newcommand{\Y}{{\boldsymbol Y}}
\newcommand{\GL}{\mathrm{GL}}

\newcommand{\Q}{\overline{\mathbb Q}}
\newcommand{\QQ}{\mathbb Q}
\newcommand{\C}{\mathbb{C}}
\newcommand{\Z}{\mathbb{Z}}

\newcommand{\e}{{\boldsymbol e}}

\newcommand{\f}{{\boldsymbol f}}
\newcommand{\g}{{\boldsymbol g}}
\newcommand{\bw}{{\boldsymbol \xi}}

\newcommand{\I}{\mathcal{I}}

\newcommand{\val}{\operatorname{val}}

\newcommand{\house}[1]{%
	\mathchoice%
	{\maisonaux\displaystyle{#1}}%
	{\maisonaux\textstyle{#1}}%
	{\maisonaux\scriptstyle{#1}}%
	{\maisonaux\scriptscriptstyle{#1}}%
}
\newcommand{\maisonaux}[2]{%
	\sbox0{$#1#2$}%
	\dimen0=\ht0%
	\advance\dimen0 by 2.5pt
	\dimen2=\dp0%
	\advance\dimen2 by 0.5pt
	\dimen4=\dimen0%
	\advance\dimen4 by -0.4pt%
	\dimen6=\wd0%
	\advance\dimen6 by 6pt
	\mathopen{\vrule height\dimen0 depth\dimen2 width 0.4pt}%
	\kern1pt%
	\rlap{\kern-1pt\raise\dimen4\hbox{\rule{\dimen6}{0.4pt}}}%
	\mkern3mu#1#2\mkern3mu%
	\kern1pt%
	\mathclose{\vrule height\dimen0 depth\dimen2 width 0.4pt}%
}

\title{The Siegel--Shidlovskii method revisited}

\author{Boris Adamczewski}
\address{
	Univ Lyon, Universit\'e Claude Bernard Lyon 1\\
	CNRS UMR 5208, Institut Camille Jordan, F-69622 Villeurbanne Cedex, France, and
	Centre International de Rencontres Math\'ematiques (CIRM), CNRS UAR 822, 163 Av.\ de Luminy, 13009
	Marseille, France}
\email{boris.adamczewski@math.cnrs.fr}

\author{Colin Faverjon}
\address{Universit\'e Picardie Jules Verne 
	\\ CNRS UMR 7352, Laboratoire Ami\'enois de Math\'ematique Fondamentale et Appliqu\'ee 
	\\ 33 rue Saint-Leu, 80039 Amiens, Fran\-ce}
\email{colin.faverjon@math.cnrs.fr}
\date{\today}

\date{}
\thanks{Funded by the European Union (ERC Advanced Grant M\^{E}TIS, Grant Agreement No. 101264813). Views and opinions expressed are however those of the authors only and do not necessarily reflect those of the European Union or the European Research Council Executive Agency. Neither the European Union nor the granting authority can be held responsible for them.}
\begin{document}

\begin{abstract}
We revisit the classical Siegel--Shidlovskii method in transcendental number theory from a new perspective, in which functional relations play a central role. Within a single self-contained framework, we recover the foundational results in the theory of $E$-functions, including the Siegel--Shidlovskii theorem, Beukers' lifting theorem, and Andr\'e's theorem on the singularities of the minimal differential operator annihilating an $E$-function. Notably, the latter is obtained without any recourse to the theory of $G$-functions. This unified approach also allows us to resolve a problem posed by Shidlovskii: we obtain a general algebraic independence measure extending the classical Siegel--Lang measure to an arbitrary algebraically independent family of $E$-values. We also provide a short elementary proof of a second theorem of Andr\'e, which has become an important tool in the theory of $E$-functions and concerns the primality of the specialization, at a regular point, of the ideal of functional relations of a vector of $E$-functions satisfying a differential system.
\end{abstract}

\maketitle

\section{Introduction}\label{sec:intro}

In his landmark memoir of 1929, Siegel~\cite{Si29} introduced a method
for proving the transcendence and algebraic independence of values at
algebraic points of certain power series solutions of linear
differential equations, which he called $E$-functions because of their
similarity with the exponential function. His aim was to extend the theory of Hermite, which culminates in
the Lindemann--Weierstrass theorem, beyond the exponential function,
a goal he achieved for the Bessel function.

A decisive breakthrough was later made by Shidlovskii, who removed the
normality assumption introduced by Siegel in his 1949
monograph~\cite{Si_Liv}, thereby eliminating the main obstacle to the
general application of the method. This led to what is now known as the
Siegel--Shidlovskii Theorem: if $f_1,\ldots,f_m$ are $E$-functions
satisfying a differential system
\begin{equation}\label{eq:diff_sys}
Y'(z)=A(z)Y(z),\qquad A\in\mathcal M_m(\Q(z)),
\end{equation}
and if $\alpha\in\Q^\times$ is regular with respect to
\eqref{eq:diff_sys}, \textit{i.e.} not a pole of $A(z)$, then
\[
{\rm tr.deg}_{\Q}\bigl(f_1(\alpha),\ldots,f_m(\alpha)\bigr)
=
{\rm tr.deg}_{\Q(z)}\bigl(f_1(z),\ldots,f_m(z)\bigr).
\]
In particular, the algebraic independence of $f_1,\ldots,f_m$ over
$\Q(z)$ implies that of $f_1(\alpha),\ldots,f_m(\alpha)$ over $\Q$.

The subsequent development of the theory followed, for the most part,
three major directions.

\begin{itemize}

\item[{\rm (i)}] The functional algebraic independence problem underlying
the application of the Siegel--Shidlovskii Theorem: developing methods
to prove the algebraic independence of solutions of differential systems
of the form~\eqref{eq:diff_sys} and applying them to concrete examples.
This led to numerous criteria, developed in particular by Shidlovskii
and his school (see \cite{Sh_Liv} and the references therein), and later to methods based on differential Galois theory,
culminating in the work of Beukers, Brownawell, and
Heckman~\cite{BBH}.

\item[{\rm (ii)}] Refinements about linear independence.
Lang conjectured~\cite[p.~100]{La_Liv} that, under the assumptions of the
Siegel--Shidlovskii Theorem, the linear independence of 
$f_1,\ldots,f_m$ over $\Q(z)$ should imply that of the numbers 
$f_1(\alpha),\ldots,f_m(\alpha)$ over $\Q$. After partial results of Shidlovskii \cite{Sh94,Sh95,Sh96,Sh99} 
and of Nesterenko and Shidlovskii~\cite{NS96}, this was finally
settled by Beukers~\cite{Be06}, who proved the more precise lifting theorem:
every homogeneous algebraic relation among
$f_1(\alpha),\ldots,f_m(\alpha)$ is obtained by specializing a
homogeneous algebraic relation over $\Q[z]$ among
$f_1,\ldots,f_m$.

\item[{\rm (iii)}] The quantitative theory, namely the search for
algebraic independence measures for values of $E$-functions. This direction goes back to
Siegel~\cite{Si29}. The first general measure, under the assumption that 
$f_1,\ldots,f_m$ are algebraically independent, was obtained by Lang in
1962~\cite{La62}. It was subsequently refined by various authors.
Shidlovskii also pursued the more intricate problem of obtaining such measures for an
algebraically independent subfamily of
$\{f_1(\alpha),\ldots,f_m(\alpha)\}$, without assuming the algebraic independence of 
$f_1,\ldots,f_m$.
See Section~\ref{sec:problem(iii)} for further details.

\end{itemize}

At the end of the 1990s, Andr\'e~\cite{An1,An2} introduced a new
perspective based on the links between $E$- and $G$-functions revealed
by the Borel and Laplace transforms within the theory of arithmetic
Gevrey series. For our purposes, the main consequence of his work is the
following; see \cite[Th.\ 4.3,~i)]{An1}.

\begin{itemize}
\item[{\rm (A1)}] \emph{The only possible non-trivial singularities of
the minimal differential operator annihilating an $E$-function are $0$
and $\infty$.
}
\end{itemize}

One of the most striking aspects of Andr\'e's approach is that,
once~{\rm (A1)} is established, the Siegel--Shidlovskii Theorem follows
from a comparatively short and essentially qualitative argument. This
contrasts sharply with the classical proofs in transcendence theory,
which rely on delicate auxiliary constructions and highly quantitative
estimates. The quantitative core of the argument has not disappeared,
however, but has been transferred to the proof of~{\rm (A1)}, which
rests on the theory of $G$-functions and, in particular, on a deep
Diophantine result due to the Chudnovskys
\cite[Thm.~VIII.1.5]{DGS}: the minimal differential equation satisfied
by a $G$-function satisfies the so-called Galochkin condition. One then shows that
any differential equation satisfying Galochkin's condition is regular
singular at $0$. This can be deduced, as presented in \cite{DGS}, by
combining results of Andr\'e and Katz, or proved more directly by an
argument of Beukers 
(see \cite[Prop.\ 6.1]{LePetit}). Using the Laplace
transform, Andr\'e~\cite{An1} was then able to deduce~{\rm (A1)} from
these properties of $G$-operators. Until the present work, this was the
only known strategy for proving~{\rm (A1)}.

The importance of~{\rm (A1)} was further confirmed by Beukers' proof of
the lifting theorem~\cite{Be06}, which relies crucially on this result,
thereby settling the problem described in Direction~{\rm (ii)}. More
generally,~{\rm (A1)} makes it possible to derive statements concerning
\emph{arbitrary values} of $E$-functions at algebraic points, independently of
the singularities of a differential system in which the functions are
embedded; see Section~\ref{sec:problem(ii)} for further details. It also
plays an essential role in Delaygue's result~\cite{De25} on the
Lindemann--Weierstrass theorem, and in the recent work of the first author
with Bell and Delaygue~\cite{ABD1,ABD2}, where a new method is developed
to prove the algebraic independence of a subfamily of
$\{f_1(\alpha),\ldots,f_m(\alpha)\}$, without first establishing the
algebraic independence over $\Q(z)$ of the whole family
$f_1,\ldots,f_m$, as in Direction~{\rm (i)}.

Another important result in the theory of $E$-functions is the following
primality statement. It follows directly from a theorem of
Andr\'e~\cite[Thm.~1.6.1(3)]{An14}, although his result is formulated in
the language of algebraic geometry.

\begin{itemize}
\item[{\rm (A2)}] \emph{Under the assumptions of the
Siegel--Shidlovskii Theorem, let $\mathfrak I$ be the homogeneous ideal
of $\Q[z][X_1,\ldots,X_m]$ generated by all homogeneous algebraic
relations over $\Q[z]$ among $f_1,\ldots,f_m$, and let
\[
{\rm ev}_\alpha:\Q[z,X_1,\ldots,X_m]
\longrightarrow \Q[X_1,\ldots,X_m]
\]
be the evaluation morphism defined by
$
{\rm ev}_\alpha(P):=P(\alpha,X_1,\ldots,X_m)$.  
\\
Then ${\rm ev}_\alpha(\mathfrak I)$ is a prime ideal of
$\Q[X_1,\ldots,X_m]$.}
\end{itemize}

In~\cite{An14}, Andr\'e used his generalized differential Galois theory
to establish~{\rm (A2)} and thereby derive a new proof of the lifting theorem
that does not rely on~{\rm (A1)}. Indeed, Statement~{\rm (A2)}, together with the
Siegel--Shidlovskii theorem, yields the lifting theorem by elementary
commutative algebra, essentially through a comparison of Krull
dimensions. Earlier, Nesterenko and
Shidlovskii~\cite[Thm.~1]{NS96} had already obtained the lifting theorem
under the assumption~{\rm (A2)} by adapting the classical
Siegel--Shidlovskii argument. Shidlovskii also obtained several further
results conditional on~{\rm (A2)} (see \cite{Sh97,Sh99,Sh00}).

\medskip

In this article, we revisit the classical Siegel--Shidlovskii method from
a new perspective, described in Section~\ref{sec:method}. Our central
idea is to assign a direct role to functional relations in the
construction, rather than attempting to avoid them. 
A single underlying argument, with only minor variations, then yields
the Siegel--Shidlovskii Theorem, Beukers' lifting theorem, and a general
algebraic independence measure. For the Siegel--Shidlovskii Theorem
itself, this perspective brings no substantial simplification. For the
lifting theorem, however, our perspective makes both the role
of~{\rm (A2)} and its close connection with the proof of the
Siegel--Shidlovskii Theorem particularly clear. It also leads to a
substantially simpler implementation of the direct approach of
Nesterenko and Shidlovskii, avoiding a technical reduction process that
is no longer needed in our framework. 
The main gain occurs in the quantitative setting. We establish a general
algebraic independence measure (see
Theorem~\ref{thm:measure_system}), which settles the problem pursued by
Shidlovskii in Direction~{\rm (iii)}. Combined with~{\rm (A1)}, it yields
an algebraic independence measure of the same general type as those of
Siegel and Lang for any algebraically independent family of
$E$-values; see Corollary~\ref{thm:measure_sans_system}. 
A further variation of the method, still based on assigning a central
role to functional relations, yields a proof of~{\rm (A1)} entirely
within the theory of $E$-functions. 
In Appendix~\ref{AppendiceA}, we also
give a short elementary proof of~{\rm (A2)}, avoiding any appeal to
differential Galois theory, either in its classical or generalized form.

Thus, apart from Siegel's classical construction based on Pad\'e
approximants and Shidlovskii's lemma, both of which are taken
for granted here, 
the present paper recovers all the fundamental results of the
transcendental theory of $E$-functions within a single self-contained
framework.
\subsection{Shidlovskii's problem on algebraic independence measures}
\label{sec:problem(iii)}

In the first part of his 1929 memoir, Siegel applied his method to the
Bessel function
\[
J_0(z):=\sum_{n\geq 0}\frac{(-1)^n}{(n!)^2}
\left(\frac{z}{2}\right)^{2n}.
\]
His primary aim was to prove that $J_0(\alpha)$ is transcendental for
every non-zero algebraic number $\alpha$, but his main result, the
\emph{Hauptsatz} of~\cite{Si29}, was already quantitative. More
precisely, if $\alpha$ is a non-zero algebraic number of degree $\delta$, and
if $P\in\mathbb Z[X,Y]$ is a non-zero polynomial of degree at most $d$
and height at most $H$, then
\[
\left|P\bigl(J_0(\alpha),J_0'(\alpha)\bigr)\right|
>
c\,H^{-123\delta^3d^2},
\]
where $c>0$ depends only on $\alpha$ and $d$. In particular,
$J_0(\alpha)$ and $J_0'(\alpha)$ are algebraically independent over
$\Q$.

The first general algebraic independence measure in the setting of a
differential system was obtained by Lang in 1962~\cite{La62}. Under the
assumptions of the Siegel--Shidlovskii Theorem, and assuming moreover
that $f_1,\ldots,f_m$ are algebraically independent over $\Q(z)$, Lang
proved that
\begin{equation}\label{eq:lang}
\left|P\bigl(f_1(\alpha),\ldots,f_m(\alpha)\bigr)\right|
\geq
C_1H^{-C_2d^m}
\end{equation}
for every non-zero polynomial
$P\in\mathbb Z[X_1,\ldots,X_m]$ of degree at most $d$
and height at most $H$. Here $C_2>0$ depends only on
$f_1,\ldots,f_m$ and $\alpha$, while, for every $d\geq 1$, the constant
$C_1>0$ depends only on $d$, $f_1,\ldots,f_m$, and $\alpha$, and may be
chosen uniformly for all polynomials of degree at most~$d$. The constants
 in this measure were subsequently refined by several
authors~\cite{Ga68,Ne77,Sh79,Sh82,Br85,Gor90}.

Shidlovskii devoted considerable effort to the more difficult problem of
studying the algebraic independence of arbitrary subfamilies
\begin{equation}\label{eq:subsystem}
\{f_{i_1},\ldots,f_{i_r}\}
\subset
\{f_1,\ldots,f_m\}
\end{equation}
and obtaining an algebraic independence measure whenever the
corresponding values are algebraically independent; see
\cite[Chaps.~11--13]{Sh_Liv} and
\cite{Sh80a,Sh80b,Sh82,Sh83}. In his last published paper~\cite{Sh00}, he
considered the still more general case in which the values need not be
algebraically independent, and formulated the desired estimate in its
natural alternative form: a polynomial in
$f_1(\alpha),\ldots,f_m(\alpha)$ either vanishes or its absolute value satisfies 
a lower bound of the form~\eqref{eq:lang}. He obtained only a partial result in this direction.

Our first main result settles Shidlovskii's problem.

\begin{thm}\label{thm:measure_system}
Under the assumptions of the Siegel--Shidlovskii Theorem and with the
notation of~\eqref{eq:subsystem}, the following alternative holds for every
polynomial $P\in\mathbb Z[X_1,\ldots,X_r]$:
\begin{equation}\label{eq:mesE}
\begin{aligned}
&\text{either}\quad
P\bigl(f_{i_1}(\alpha),\ldots,f_{i_r}(\alpha)\bigr)=0,\\
&\text{or its absolute value satisfies a lower bound of the
form~\eqref{eq:lang},}
\end{aligned}
\end{equation}
with the exponent $m$ replaced by
$
t:={\rm tr.deg}_{\Q(z)}
\bigl(f_1(z),\ldots,f_m(z)\bigr)\leq m$.
\end{thm}

The vanishing alternative is essential since, in this generality, the
corresponding values may well be algebraically dependent. When
$f_{i_1}(\alpha),\ldots,f_{i_r}(\alpha)$ are algebraically independent,
however, 
Theorem~\ref{thm:measure_system} becomes an algebraic independence
measure in the usual sense.


\subsection{Andr\'e's theorem and arbitrary values of \(E\)-functions}
\label{sec:problem(ii)}

Beyond its crucial role in Beukers' proof of the lifting
theorem~\cite{Be06}, Theorem~{\rm (A1)} has been the starting point of
many recent developments concerning arbitrary values of
$E$-functions at algebraic points. Let
\[
\mathbf E
:=
\{f(\alpha): f\text{ is an $E$-function},\ \alpha\in\Q\}
\]
denote the ring of values of $E$-functions at algebraic points and, for
a subfield $\mathbb K\subset\Q$, let $\mathbf E_{\mathbb K}$ denote the
set of values $f(\alpha)$ with $f\in\mathbb K[[z]]$ and
$\alpha\in\mathbb K$. Any statement applying uniformly to all elements
of $\mathbf E$, or of $\mathbf E_{\mathbb K}$, necessarily requires a
means of overcoming the possible singularities of the differential
systems involved; this is precisely what~{\rm (A1)} provides. We list below some of the results made possible 
by this theorem, notably in the work of Fischler, Rivoal, and the authors.

\begin{itemize}

\item Elements of $\mathbf E_{\mathbb K}$ are linearly dependent over
$\Q$ if and only if they are linearly dependent over $\mathbb K$. More
precisely, if $\mathbb K\subset\mathbb L\subset\Q$ are fields, then the canonical map $\mathbf E_{\mathbb K}\otimes_{\mathbb K}\mathbb L \to \mathbf E_{\mathbb L}$, $\xi \otimes \alpha \mapsto \alpha\xi$, is an isomorphism~\cite{AF24,FR16,FR23}.

\item There is an algorithm to determine all the algebraic values taken
by a given $E$-function at algebraic points. More generally, one can
decide linear dependence among elements of $\mathbf E$ and, modulo the
computation of the corresponding functional ideal, determine all their
algebraic relations~\cite{AF24,AR18,BRS,FR19}.

\item Every algebraic relation among elements of $\mathbf E$ has a
functional origin: it is obtained by specializing an algebraic
differential relation among the functions involved~\cite{AF24}.

\item Non-zero linear forms in elements of $\mathbf E$ satisfy
quantitative lower bounds. In particular, no element of $\mathbf E$ is
a Liouville number, and every irrational element of
$\mathbf E_{\mathbb Q}$ has irrationality exponent equal to
$2$~\cite{FR23,FR23_Exp}.

\end{itemize}

Combining Andr\'e's theorem~{\rm (A1)} with
Theorem~\ref{thm:measure_system}, we obtain the following intrinsic
result.

\begin{coro}\label{thm:measure_sans_system}
Let $\xi_1,\ldots,\xi_r\in\mathbf E$ be algebraically independent. Then,
for every non-zero polynomial $P\in\mathbb Z[X_1,\ldots,X_r]$, the
quantity
$
\left|P(\xi_1,\ldots,\xi_r)\right|$ 
satisfies a lower bound of the form~\eqref{eq:lang}, with the exponent
$m$ replaced by some positive integer $t$ independent of $P$.
\end{coro}

Thus, as soon as values of $E$-functions are known to be algebraically
independent, they \emph{automatically} satisfy an algebraic independence
measure of Siegel--Lang type. As shown in
Section~\ref{sec:measures_E}, the exponent $t$ may moreover be bounded
in terms of the transcendence degree of the differential field
generated by the $E$-functions used to represent
$\xi_1,\ldots,\xi_r$. The fact that no element of $\mathbf E$ is a
Liouville number also follows immediately from
Corollary~\ref{thm:measure_sans_system}.

\subsection{Organization of the paper}

In Section~\ref{sec:method}, we recall the classical
Siegel--Shidlovskii method and present our new perspective.
Section~\ref{sec:systems} deals with values of $E$-functions at regular
points of differential systems. There we give direct proofs of the
Siegel--Shidlovskii Theorem and Beukers' lifting theorem, and establish
the general algebraic independence measure of
Theorem~\ref{thm:measure_system}. In Section~\ref{sec:Andre}, we extend
the method to points that are singular for the differential system
under consideration and obtain a new proof of Andr\'e's
theorem~{\rm (A1)}. We then deduce
Corollary~\ref{thm:measure_sans_system}.
Appendix~\ref{AppendiceA} is devoted to a proof of Andr\'e's
result~{\rm (A2)}. Finally, in Appendix~\ref{AppendiceB}, we discuss
consequences of our approach for the local structure of minimal
differential equations annihilating $E$- and $G$-functions.

\subsection{Notation and definitions}

We let $\Q$ denote the algebraic closure of $\mathbb Q$ in $\C$. 
Given a number field $\mathbb K \subset \Q$, we let $\OK$ denote its ring of integers. 
The house of an algebraic number $\alpha$ is the maximum modulus of its
Galois conjugates and is denoted by $\house{\alpha}$. The house of a linear form $L$ with algebraic coefficients, denoted 
by $\house{L}$, is the maximum of the houses of its coefficients. 
We recall that a positive integer $d$ is called a denominator of an
algebraic number $\alpha$ if $d\alpha$ is an algebraic integer. 
As usual in the theory of $E$-functions, and throughout this paper, denominators and houses are used to measure the arithmetic complexity of algebraic numbers and linear forms, playing the role more commonly assigned to height functions. 
To align with standard notation, given a polynomial $P \in \Z[X_1,\ldots,X_m]$, we let $H(P)$ denote its height, that is, the maximum of the absolute values of its coefficients. In this setting, $H(P)$ coincides with $\house{P}$ when $P$ is viewed 
as a linear form in its monomials.

Given a linear differential equation
\[
a_m(z)Y^{(m)}+a_{m-1}(z)Y^{(m-1)}+\cdots+a_0(z)Y=0,
\]
with coefficients in $\mathbb C[z]$, a point $\alpha\in\mathbb C$ is
\emph{singular} if it is a zero of the polynomial $a_m(z)$. A point that is
not singular is called a \emph{regular} or an \emph{ordinary} point. A \emph{trivial} or \emph{apparent} singularity is a singular point at which the
equation still admits a local basis of holomorphic solutions.
A singular point $\alpha$ is \emph{regular singular} if, for every
$i\in\{1,\ldots,m\}$, the rational function $(z-\alpha)^{i}a_{m-i}(z)/a_m(z)$
is holomorphic at $\alpha$. A singularity which is not regular singular is called \emph{irregular}.

In the case of a linear system
\[
Y'=A(z)Y,\qquad A(z)\in\mathcal M_m(\mathbb C(z)),
\]
a point $\alpha\in\mathbb C$ is singular if it is a pole of at least one of the entries of
$A(z)$, and regular or ordinary otherwise.

Finally, recall that a power series 
$$
 \sum_{n=0}^\infty \frac{a_n}{n!} z^n \in \Q[[z]]
$$ is an $E$-function if it satisfies a linear differential equation over $\Q(z)$ and if there exists $C>0$ such that $\house{a_n} \leq C^{n}$ and $d_n \leq C^{n}$ for all positive integers $n$, where $d_n$ is the least common denominator of $a_0,\ldots,a_n$. The field generated by the coefficients $a_n$, $n\geq 0$, of an $E$-function is a number field.

\section{The Siegel--Shidlovskii method}\label{sec:method}

The purpose of this section is to recall the main steps of the Siegel--Shidlovskii 
method and then to present our new perspective.

\subsection{Siegel's Diophantine argument}

The first objective of the method is to prove the linear independence of
complex numbers, say $\xi_1,\ldots,\xi_\nu$, or at least to obtain a lower bound for their
linear rank. By viewing algebraic relations as linear forms in monomials, this framework
naturally extends to algebraic independence and lower bounds on the
transcendence degree. Natural refinements also concern linear and
algebraic independence measures.

In each case, one argues by contradiction. The failure of the desired
conclusion gives rise, for some positive integer $\eta$,  to linear forms
$
L_1,\ldots,L_\eta$ 
with coefficients in $\mathcal O_{\mathbb K}$ for some number field $\mathbb K$, linearly independent over $\mathbb K$,
and taking \emph{zero or exceptionally small values} at
\[
\bw:=\prescript{\rm t}{}{(\xi_1,\ldots,\xi_\nu)}
\in\mathbb C^\nu\setminus\{0\}.
\]
We call $L_1,\ldots,L_\eta$ the \emph{putative forms}.

Independently of these forms, the method constructs, for large values of
an auxiliary parameter $n$, a basis of linear forms in $\nu$ variables
with coefficients of controlled arithmetic complexity and taking small
values at $\bw$.  
Their coefficient vectors form an invertible $\nu \times \nu$
matrix
$
M_n$ with coefficients in $\mathcal O_{\mathbb K}$.
By construction, 
$
M_n\bw=\mathbf e_n$,
where $\mathbf e_n$ is a vector with \emph{small norm}.   
We call this matrix a \emph{preparation matrix}.

One then replaces suitable rows of $M_n$ by the coefficient vectors of
the putative forms while preserving invertibility. The
resulting matrix
$
M'_n\in\mathcal M_\nu(\mathcal O_{\mathbb K})
$
satisfies
$
M'_n\bw=\boldsymbol\varepsilon_n$,
where some coordinates of $\boldsymbol\varepsilon_n$ are much smaller than the corresponding ones in $\mathbf e_n$.
Since $M'_n$ is invertible, one has
\begin{equation}\label{eq:inverse}
\bw=(M'_n)^{-1}\boldsymbol\varepsilon_n.
\end{equation}

By the comatrix formula, we have
\[
(M'_n)^{-1}
=
\frac{1}{\det M'_n}\prescript{\rm t}{}{\operatorname{Com}(M'_n)},
\]
where the entries of $\operatorname{Com}(M'_n)$ are signed
$(\nu-1)\times(\nu-1)$ minors of $M'_n$. The upper bounds for the entries
of $M'_n$ therefore control these minors. Moreover,
$
\det M'_n$ is a non-zero algebraic integer, 
so Liouville's inequality gives a lower bound for $|\det M'_n|$ in terms
of its house. Combining these estimates yields, for each $j\in\{1,\ldots,\nu\}$, a bound
\[
\max_{1\leq i\leq \nu}
\left|((M'_n)^{-1})_{ij}\right|
\leq C_{j,n},
\]
where the quantity $C_{j,n}$ may depend on the column. Writing
$$\boldsymbol\varepsilon_n:={}^t(\varepsilon_{1,n},\ldots,\varepsilon_{\nu,n}),$$
one obtains
\[
\|(M'_n)^{-1}\boldsymbol\varepsilon_n\|_\infty
\leq
\sum_{j=1}^{\nu} C_{j,n}|\varepsilon_{j,n}|.
\]
Since $\bw\neq0$, if, for a suitable choice of $n$, the values of the
putative forms make the right-hand side sufficiently small, one obtains
\[
\|(M'_n)^{-1}\boldsymbol\varepsilon_n\|_\infty
<
\|\bw\|_\infty,
\]
in contradiction with~\eqref{eq:inverse}. We call a matrix $M'_n$ for 
which this contradiction occurs a \emph{forbidden matrix}.

The method can therefore be summarized in three steps. One first
constructs a sequence of preparation matrices $(M_n)_n$, then replaces
some of their rows by coefficient vectors of putative forms, and finally
proves that, for a suitable $n$, the resulting matrix $M'_n$ is
a forbidden matrix. The essential difficulty lies in the first step.

\begin{rem*}
The preparation matrices $M_n$ need not be invertible; it is enough that
the matrices $M_n'$ are. We shall use this when proving Andr\'e's theorem~{\rm (A1)}.
\end{rem*}

\subsection{Classical construction of the preparation matrix}

Let $f_1,\ldots,f_m$ be $E$-functions which are linearly independent over
$\overline{\mathbb Q}(z)$ and satisfy
\begin{equation*}
Y'(z)=A(z)Y(z),
\qquad
A(z)\in\mathcal M_m(\overline{\mathbb Q}(z)).
\end{equation*}
Let $\alpha\in\overline{\mathbb Q}^{\times}$ be a regular point of this
system, and set
\[
\bw:=\prescript{\rm t}{}{(f_1(\alpha),\ldots,f_m(\alpha))}.
\]
We fix a number field $\mathbb K$ containing $\alpha$, the coefficients
of $A(z)$, and all the algebraic coefficients occurring in the
construction.

The classical preparation matrix is obtained by combining Siegel's
lemma and Shidlovskii's lemma. Siegel's lemma, applied over $\mathbb K$,
provides polynomials
$
P_{1,n}(z),\ldots,P_{m,n}(z)\in\OK[z]$ 
of controlled degree, denominator, and house such that the functional linear form 
\[
R_n(z)
:=
P_{1,n}(z)f_1(z)+\cdots+P_{m,n}(z)f_m(z)
\]
has a large order of vanishing at the origin.

The differential system allows one to differentiate $R_n$ and, after clearing the denominators, to
rewrite its successive derivatives in the form
\[
R_{n,k}(z)
:=
P_{1,n,k}(z)f_1(z)+\cdots+P_{m,n,k}(z)f_m(z),
\]
where the polynomials $P_{i,n,k}$ have coefficients in $\OK$ and controlled degree and house.

Shidlovskii's lemma ensures that, for a second auxiliary parameter
$\varepsilon>0$, which may be chosen arbitrarily small, one can select
$m$ such forms among those with indices
$1\leq k\leq \varepsilon n+C$ whose specialized coefficient vectors at
$z=\alpha$ are linearly independent. Here $C$ is independent of $n$.
This bound on the indices guarantees that their coefficients remain
controlled and that their values at $\bw$ remain small. 
These
vectors form the rows of an invertible preparation matrix
$M_n\in\mathcal M_m(\mathbb K)$.

Two fundamental restrictions are inherent in this construction:
\begin{enumerate}
\item the functions must be linearly independent over
$\overline{\mathbb Q}(z)$;

\item the point $\alpha$ must be regular for the differential system.
\end{enumerate}

When $\mathbb K = \mathbb Q$ or $\mathbb K$ is an imaginary quadratic field, the computation shows that one single putative form is sufficient to obtain a forbidden matrix. The linear independence of $f_1(\alpha),\ldots,f_m(\alpha)$ over 
$\mathbb K$ follows. 
If the functions $f_1,\ldots,f_m$ are algebraically independent over $\Q(z)$, one deduces the algebraic independence of $f_1(\alpha),\ldots,f_m(\alpha)$
over $\Q$ by applying the same argument to the associated monomials and
using their linear independence.  For general number fields, however, more than one putative form is needed, and the contradiction only yields a lower bound for the $\mathbb K$-rank of $f_1(\alpha),\ldots,f_m(\alpha)$. But one can still recover the algebraic independence of $f_1(\alpha),\ldots,f_m(\alpha)$  by applying such a lower bound asymptotically to monomials of large degree.

\subsubsection{Obstruction arising from functional relations}

When $f_1,\ldots,f_m$ are linearly dependent, the classical remedy is to choose a basis of their
span and to apply the method to the induced differential system, thus
reducing the dimension. 
This reduction changes the natural coordinates of the problem. More
precisely, the preparation matrix is no longer constructed from the
original family of functions, but from a basis of its
$\overline{\mathbb Q}(z)$-span. Thus, if the initial family consists of
$m$ functions of rank $r<m$, the preparation matrix has size
$r\times r$ rather than $m\times m$. 

The same phenomenon becomes much
more pronounced when studying algebraic relations. There, one starts
with the monomials of degree $\delta$ in the functions
$f_1,\ldots,f_m$, but the preparation matrix is again constructed only
after replacing this natural family by a basis of its
$\overline{\mathbb Q}(z)$-span. The putative forms, however, are
naturally expressed in terms of the original monomials. One must
therefore rewrite them in the chosen basis and control their
independence throughout this change of coordinates. This reduction
process is a major source of technical complications and has proved to
be a serious obstacle in several situations.

\subsection{The new construction}

The key difference with the classical method lies in the treatment of
functional relations. Instead of eliminating them before constructing the
preparation matrix, our approach incorporates them into the preparation
matrix itself. 

\subsubsection{The construction at regular points}

Let
$
g_1,\ldots,g_\nu\in\mathbb K[[z]]
$ be $E$-functions which satisfy
\begin{equation*}
Y'(z)=B(z)Y(z),
\qquad
B(z)\in\mathcal M_\nu(\overline{\mathbb Q}(z)).
\end{equation*}
Let $\alpha\in\overline{\mathbb Q}^{\times}$ be a regular point of this
system. 
One can think of $g_1,\ldots,g_\nu$ as monomials in some algebraically dependent $E$-functions $f_1,\ldots,f_m$ one 
starts with, so that $g_1,\ldots,g_\nu$ may not be linearly independent. 

Let $\mu$ denote the rank of $g_1,\ldots,g_\nu$ over $\mathbb K(z)$. The space of functional linear
relations among them has dimension $\nu-\mu$ over $\mathbb K(z)$. We
choose a basis of this space whose coefficients are regular at
$z=\alpha$ and whose specializations at $\alpha$ remain linearly
independent and have coefficients in $\OK$. We show that the Pad\'e construction, applied to the $g_i$, produces $\mu$
independent forms in directions transverse to these relations, and we
complete them by adjoining the relations themselves.

After specialization at $z=\alpha$, the preparation matrix $M_n$ has the form
\[
\begin{pmatrix}
\text{\(\mu\) Pad\'e-type forms depending on \(n\)}\\
\hline
\text{\(\nu-\mu\) specialized functional relations}
\end{pmatrix}
\in\mathcal M_\nu(\OK),
\]
and
\[
M_n\bw
=
\prescript{\rm t}{}{(e_{1,n},\ldots,e_{\mu,n},0,\ldots,0)},
\qquad
\bw:=\prescript{\rm t}{}{(g_1(\alpha),\ldots,g_\nu(\alpha))}.
\]

The functional relations supply exactly the rows missing from the
Pad\'e construction. Their coefficients are independent of $n$, and
their values at $\bw$ vanish identically. In this sense, the forms
arising from functional relations are far better than those produced
by Pad\'e approximation and differentiation. What constitutes an obstruction
in the classical construction becomes a source of exceptionally good
rows in the new one.

In the proofs of our main results, we show that the 
negation of the desired conclusion provides
 putative forms whose coefficient vectors  can be substituted for  
Pad\'e-type rows. Since all the forms remain expressed in
the original coordinates, their independence is preserved, and the
resulting matrix is a forbidden matrix.

This principle is the key to the proof of
Theorem~\ref{thm:measure_system}. It also clarifies and unifies the
arguments: Theorem~\ref{thm:measure_system}, the
Siegel--Shidlovskii Theorem, and Beukers' lifting theorem all follow
from essentially the same mechanism. The preparation matrices and the
forbidden-matrix argument remain unchanged; only the putative forms and
the final conclusion differ. 


\subsubsection{The construction at singular points}

The same principle makes it possible to work at singular points.

At a regular point, the successive derivatives of a Pad\'e approximant can
be rewritten as linear forms in the original functions, with
coefficients regular at the point of specialization. At a singular
point, the differential system introduces poles and therefore no longer
directly produces a preparation matrix.

The key idea is therefore to avoid this rewriting procedure by simply
enlarging the family of functions, adjoining the derivatives that
occur in the Pad\'e construction. The resulting family is linearly
dependent over $\mathbb K(z)$ and would therefore fall outside the scope
of the classical method. In our approach, however, the functional
relations simply provide additional rows of a higher-dimensional
preparation matrix.

Suppose that an $E$-function $f$ satisfies a linear differential
equation of order $\ell$ which is singular at $\alpha$. Then the
corresponding relation among $f,f',\ldots,f^{(\ell)}$ specializes at
$z=\alpha$ to a non-trivial relation of smaller order, involving only
\[
f(\alpha),\ldots,f^{(\ell-1)}(\alpha).
\]
In particular, if $f$ satisfies no linear differential equation which
is regular at $\alpha$, then, after adjoining sufficiently many
derivatives, one obtains too many specialized relations of smaller order
to be accommodated by the preparation matrix. Inserting these
relations as rows ultimately produces a forbidden matrix. This yields
our proof of Andr\'e's theorem~{\rm (A1)}.

\section{Applying the method at regular points}
\label{sec:systems}

In this section, we consider values of $E$-functions related by a
differential system at a non-zero algebraic point which is regular for
this system. We prove the Siegel--Shidlovskii Theorem, Beukers'
lifting theorem, and Theorem~\ref{thm:measure_system}.

All three proofs rely on the same sequence of preparation matrices. 
They only differ in the construction of the
putative forms. In the proofs of the Siegel--Shidlovskii and lifting theorems, the
putative forms vanish at the point under consideration. In the proof of
the algebraic independence measure, they are non-zero but extremely small.
This requires a more delicate adjustment of the auxiliary parameters
and makes the quantitative argument somewhat longer, whereas the
underlying matrix argument remains the same.

\subsection{Notation}
\label{sec:notation_method}
Throughout this section, let $\mathbb K$ be a number field, let
$f_1,\ldots,f_m\in\mathbb K[[z]]$ be $E$-functions satisfying a
differential system
\[
Y'(z)=A(z)Y(z),
\qquad
A(z)\in\mathcal M_m(\mathbb K(z)),
\]
and let $\alpha\in\mathbb K^\times$ be a point which is not a pole of
$A$. Since the class of $E$-functions is invariant under the change of
variable $z\mapsto\alpha z$, we may assume that $\alpha=1$.

We denote by $\mathfrak I$ the ideal of
$\mathbb K[z,X_1,\ldots,X_m]$ generated by all polynomials
$Q(z,X_1,\ldots,X_m)$ that are homogeneous in
$X_1,\ldots,X_m$ and satisfy
\[
Q(z,f_1(z),\ldots,f_m(z))=0.
\]
Thus $\mathfrak I$ is a homogeneous ideal with respect to the variables
$X_1,\ldots,X_m$. We further set
\[
\mathfrak I_{\mathbb K(z)}
:=
\mathbb K(z)\mathfrak I
\subseteq
\mathbb K(z)[X_1,\ldots,X_m].
\]
Then $\mathfrak I_{\mathbb K(z)}$ is precisely the ideal generated by
the homogeneous algebraic relations over $\mathbb K(z)$ among
$f_1,\ldots,f_m$.

We recall that ${\rm ev}_1$ denotes the specialization map
\[
{\rm ev}_1:
\begin{cases}
\mathbb K[z,X_1,\ldots,X_m]\longrightarrow
\mathbb K[X_1,\ldots,X_m],\\
Q(z,X_1,\ldots,X_m)\longmapsto Q(1,X_1,\ldots,X_m).
\end{cases}
\]
Finally, let $\mathfrak I(1)$ denote the ideal of
$\mathbb K[X_1,\ldots,X_m]$ generated by all homogeneous polynomials
$P(X_1,\ldots,X_m)$ such that
\[
P(f_1(1),\ldots,f_m(1))=0.
\]
Thus $\mathfrak I(1)$ is homogeneous, and we have
$
{\rm ev}_1(\mathfrak I)
\subseteq
\mathfrak I(1)$.

Given an ideal $\mathfrak a$ homogeneous with respect to
$X_1,\ldots,X_m$ and an integer $\delta\geq0$, we denote by
$\mathfrak a_\delta$ its homogeneous component of degree $\delta$ in
these variables.

\subsection{The sequence of preparation matrices}
\label{sec:prepa_matrix}

Let $\delta\geq1$ be an integer, to be chosen later. We consider the
monomials of total degree $\delta$ in the variables
$X_1,\ldots,X_m$. Let
$Y_1,\ldots,Y_\nu$ be an enumeration of these monomials, and let
$g_1,\ldots,g_\nu$ be the corresponding monomials in
$f_1,\ldots,f_m$. We set
\[
\bw:={}^{\rm t}(g_1(1),\ldots,g_\nu(1)),
\qquad
\X:=(X_1,\ldots,X_m),
\qquad
\Y:=(Y_1,\ldots,Y_\nu),
\]
keeping in mind that each $Y_i$ is a monomial of total degree
$\delta$ in $\X$.

\subsubsection{Forms arising from functional relations}
\label{subsec:vanishing}

Let $\mu$ denote the dimension over $\mathbb K(z)$ of the vector space
spanned by $g_1,\ldots,g_\nu$. The kernel of the map
\[
\mathbb K(z)Y_1\oplus\cdots\oplus\mathbb K(z)Y_\nu
\longrightarrow
\mathbb K(z)g_1+\cdots+\mathbb K(z)g_\nu,
\qquad
Y_i\longmapsto g_i,
\]
identifies with $(\mathfrak I_{\mathbb K(z)})_\delta$. Hence
\[
\dim_{\mathbb K(z)}
(\mathfrak I_{\mathbb K(z)})_\delta=\nu-\mu.
\]
Thus $\nu-\mu$ is the number of independent homogeneous functional
relations of degree $\delta$, whereas $\mu$ is the number of
additional Pad\'e forms needed to complete the preparation matrix.

We may choose forms
\[
\ell_1(z,\Y),\ldots,\ell_{\nu-\mu}(z,\Y)
\in\mathfrak I_\delta
\]
with coefficients in $\mathcal O_{\mathbb K}[z]$ which form a basis of
$(\mathfrak I_{\mathbb K(z)})_\delta$ over $\mathbb K(z)$ and whose
specializations
\[
\ell_1(1,\Y),\ldots,\ell_{\nu-\mu}(1,\Y)
\]
are linearly independent over $\mathbb K$; see, for instance,
\cite[Lem.~3.1]{Be06}. In particular,
\[
\ell_1(1,\bw)=\cdots=\ell_{\nu-\mu}(1,\bw)=0.
\]

\subsubsection{Forms arising from Pad\'e approximants}

We now complete the family
$$\ell_1(1,\Y),\ldots,\ell_{\nu-\mu}(1,\Y)$$ by linear forms taking
small values at $\bw$. We use the following classical consequence of
Siegel's construction and Shidlovskii's lemma.

\begin{prop}
\label{prop:SiegelShidlovskii}
Let $s$ be a positive integer, and let
$h_1,\ldots,h_s\in\mathbb K[[z]]$ be linearly independent
$E$-functions related by a linear differential system of size $s$. Let
$0<\varepsilon_0<1/2$. There exists a positive real number $\rho$,
depending only on $h_1,\ldots,h_s$ and $\varepsilon_0$, such that, for
every positive integer $n$, there exist polynomials
$
P_1,\ldots,P_s\in{\mathcal O}_{\mathbb K}[z]$, 
not all zero, of degree at most $n$ and house at most
$\rho n^{(1+\varepsilon_0)n}$, such that
\[
\val R(z)\geq s(n+1)-\varepsilon_0n-1,
\]
where
\[
R(z):=P_1(z)h_1(z)+\cdots+P_s(z)h_s(z).
\]
Furthermore, there exists a positive integer $n_0$, depending only on
$h_1,\ldots,h_s$, such that
$
R(z),R'(z),\ldots,R^{(s-1)}(z)$ 
are linearly independent over $\mathbb K(z)$ whenever
$
n\geq\max\{n_0,\varepsilon_0^{-1}\}$. 
\end{prop}

\begin{proof}
The existence of $P_1,\ldots,P_s$, not all zero, satisfying
the stated degree, height, and vanishing estimates follows directly
from Siegel's lemma; see \cite[Chap.~3, Lem.~14]{Sh_Liv}. Since
$\varepsilon_0<1$, the linear independence assertion follows
immediately from Shidlovskii's lemma; see \cite[Lem.~5.2]{FN} or
\cite[Chap.~3, Lem.~8]{Sh_Liv}.
\end{proof}

The $\mathbb K(z)$-vector space spanned by $g_1,\ldots,g_\nu$ has
dimension $\mu$. After reindexing, we may assume that
$g_1,\ldots,g_\mu$ form a basis of this space. Since
$f_1,\ldots,f_m$ satisfy a linear differential system, the same holds
for $g_1,\ldots,g_\nu$, and hence for $g_1,\ldots,g_\mu$.

Let $0<\varepsilon<1$ be sufficiently small, to be chosen later. We
apply Proposition~\ref{prop:SiegelShidlovskii} with
\[
s:=\mu,
\qquad
\varepsilon_0:=\frac{\varepsilon}{6(\mu+1)},
\qquad
h_i:=g_i\quad(1\leq i\leq\mu).
\]
Set
$
n_1:=\max\{n_0,\varepsilon_0^{-1}\}$. 
For every $n\geq n_1$, let $P_1,\ldots,P_\mu$ and $R$ be provided by
Proposition~\ref{prop:SiegelShidlovskii}. 
Since $z=1$ is regular for the differential system satisfied by
$f_1,\ldots,f_m$, the functions $g_1,\ldots,g_\nu$ satisfy a
differential system
\[
Y'(z)=B(z)Y(z)
\]
which is regular at $z=1$. Let
$T(z)\in\mathcal O_{\mathbb K}[z]$ be a least common denominator of the
entries of $B(z)$, chosen so that $T(z)B(z)$ has entries in
$\mathcal O_{\mathbb K}[z]$ and $T(1)\neq0$.

Let $\partial=d/dz$. For every linear form
$L(z,\Y)\in\mathcal O_{\mathbb K}[z,\Y]$, there exist polynomials
$B_1(z),\ldots,B_\nu(z)\in\mathcal O_{\mathbb K}[z]$ such that
\[
T(z)\partial\bigl(L(z,\g(z))\bigr)
=
\sum_{i=1}^{\nu}B_i(z)g_i(z).
\]
We define
\[
\Theta(L)(z,\Y)
:=
\sum_{i=1}^{\nu}B_i(z)Y_i.
\]
Thus
$
\Theta(L)(z,\g(z))
=
T(z)\partial\bigl(L(z,\g(z))\bigr)$. 
We then define inductively
\[
R_0(z,\Y):=P_1(z)Y_1+\cdots+P_\mu(z)Y_\mu,
\qquad
R_{k+1}(z,\Y):=\Theta(R_k)(z,\Y).
\]
Notice that $R_0$ is supported on $Y_1,\ldots,Y_\mu$, whereas the
forms $R_k$, for $k\geq1$, may involve all the variables
$Y_1,\ldots,Y_\nu$.

\begin{lem}
\label{lem:formes_Pade}
There exist positive constants $c_1,\rho_1,\rho_2$, independent of
$n$, such that, for every $n\geq n_1$, there exists an integer $c_2$
satisfying
$
0\leq c_2\leq\varepsilon_0 n+c_1$ 
and $\mu$ forms
$
L_1(\Y),\ldots,L_\mu(\Y)$ 
chosen among
$
R_0(1,\Y),\ldots,R_{\mu+c_2-1}(1,\Y)$ 
such that:
\begin{itemize}
\item[\rm (a)]
$
\bigl|L_i(\bw)\bigr|
\leq \rho_1n^{-(\mu-1-\varepsilon)n}
\qquad(1\leq i\leq\mu)$;

\item[\rm (b)]
$
\house{L_i}\leq\rho_2n^{(1+\varepsilon)n}
\qquad(1\leq i\leq\mu)$;

\item[\rm (c)] the $\nu$ linear forms
$
L_1(\Y),\ldots,L_\mu(\Y),
\ell_1(1,\Y),\ldots,\ell_{\nu-\mu}(1,\Y)$ 
are linearly independent over $\mathbb K$.
\end{itemize}
\end{lem}

\begin{proof}
Once $c_1$ is fixed, the classical estimates show
that there exist positive constants $\rho_1,\rho_2$, independent of
$n$, such that every form $R_i(1,\Y)$ with
$0\leq i\leq\varepsilon_0 n+c_1$ satisfies
\[
\bigl|R_i(1,\bw)\bigr|
\leq \rho_1n^{-(\mu-1-\varepsilon)n}
\quad\mbox{ and }\quad 
\house{R_i(1,\Y)}
\leq\rho_2n^{(1+\varepsilon)n}.
\]
Detailed estimates may be found in
\cite[Chap.~3, \S10]{Sh_Liv}. Thus any forms chosen in the indicated
range satisfy~{\rm (a)} and~{\rm (b)}.

To prove~{\rm (c)}, it is enough to find an integer
$c_2\leq\varepsilon_0 n+c_1$ such that the forms
\[
R_0(1,\Y),\ldots,R_{\mu+c_2-1}(1,\Y),
\ell_1(1,\Y),\ldots,\ell_{\nu-\mu}(1,\Y)
\]
have rank $\nu$. We adapt Siegel's original argument
\cite[p.~219]{FN} to the presence of these two families of forms.

Consider the matrix
\[
\Gamma(z):=(\gamma_{i,j}(z))_{1\leq i,j\leq\nu}
\]
defined by
\[
\gamma_{i,j}(z)
:=
\begin{cases}
\text{the coefficient of $Y_j$ in $R_{i-1}(z,\Y)$},
 & 1\leq i\leq\mu,\\
\text{the coefficient of $Y_j$ in $\ell_{i-\mu}(z,\Y)$},
 & \mu<i\leq\nu.
\end{cases}
\]
It has entries in $\mathcal O_{\mathbb K}[z]$ and satisfies
\begin{equation}
\label{eq:Gamma}
\Gamma(z)
\begin{pmatrix}
g_1(z)\\
\vdots\\
g_\nu(z)
\end{pmatrix}
=
\begin{pmatrix}
R_0(z,\g(z))\\
\vdots\\
R_{\mu-1}(z,\g(z))\\
0\\
\vdots\\
0
\end{pmatrix},
\end{equation}
as well as
\[
\Gamma(z)
\begin{pmatrix}
Y_1\\
\vdots\\
Y_\nu
\end{pmatrix}
=
\begin{pmatrix}
R_0(z,\Y)\\
\vdots\\
R_{\mu-1}(z,\Y)\\
\ell_1(z,\Y)\\
\vdots\\
\ell_{\nu-\mu}(z,\Y)
\end{pmatrix}.
\] 
Let us first prove that $\Gamma(z)$ is invertible. By construction, the
functions
\[
R_0(z,\g(z)),\ldots,R_{\mu-1}(z,\g(z))
\]
span the same $\mathbb K(z)$-vector space as
\[
R(z),R'(z),\ldots,R^{(\mu-1)}(z).
\]
They are therefore linearly independent by
Proposition~\ref{prop:SiegelShidlovskii}. Moreover, the forms
\[
\ell_1(z,\Y),\ldots,\ell_{\nu-\mu}(z,\Y)
\]
are linearly independent and vanish at $\Y=\g(z)$. It follows that
\[
R_0(z,\Y),\ldots,R_{\mu-1}(z,\Y),
\ell_1(z,\Y),\ldots,\ell_{\nu-\mu}(z,\Y)
\]
are linearly independent over $\mathbb K(z)$. Hence $\Gamma(z)$ is
invertible.

Set
$
\Delta(z):=\det\Gamma(z)\neq0$. 
Since $P_1,\ldots,P_\mu$ have degree at most $n$, whereas the degrees
of the coefficients of the forms $\ell_i$ are bounded independently
of $n$, there exists a constant $\chi_1$, independent of $n$, such
that
$
\deg\Delta\leq\mu n+\chi_1$. 

On the other hand, by applying the comatrix formula to
\eqref{eq:Gamma} and using
\[
\val\bigl(R_k(z,\g(z))\bigr)
\geq\mu(n+1)-\varepsilon_0 n-k-1,
\]
we find an index $i\in\{1,\ldots,\nu\}$ such that
\[
\val(\Delta)+\val(g_i)
\geq\mu n-\varepsilon_0 n.
\]
Since the quantities $\val(g_i)$ are bounded independently of $n$, it
follows that the order of vanishing of $\Delta$ at $z=1$ is at most
$\varepsilon_0 n+\chi_2$ for some constant $\chi_2$ independent of $n$.

Let $\chi_3$ be the smallest non-negative integer such that
$
\partial^{\chi_3}\Delta(1)\neq0$. 
Then
\begin{equation}
\label{eq:rho34}
\chi_3\leq\varepsilon_0 n+\chi_2.
\end{equation}
Using the comatrix formula for $\Gamma$, we obtain, for every
$i\in\{1,\ldots,\nu\}$,
\[
\Delta(z)Y_i
=
\sum_{j=0}^{\mu-1}\Delta_{i,j+1}(z)R_j(z,\Y)
+
\sum_{j=1}^{\nu-\mu}
\Delta_{i,j+\mu}(z)\ell_j(z,\Y),
\]
where $\Delta_{i,j}(z)$ denotes the corresponding cofactor. Applying
$\Theta$ successively $\chi_3$ times and then specializing at
$z=1$, all the terms involving $\partial^k\Delta(1)$ with
$k<\chi_3$ vanish. Since $T(1)\neq0$ and
$\partial^{\chi_3}\Delta(1)\neq0$, it follows that every $Y_i$ belongs
to the $\mathbb K$-vector space spanned by
\[
R_j(1,\Y)
\qquad
(0\leq j\leq\mu+\chi_3-1)
\]
and
\[
\ell_j^{[k]}(1,\Y)
\qquad
(1\leq j\leq\nu-\mu,\ 0\leq k\leq\chi_3),
\]
where
$
\ell_j^{[0]}=\ell_j$ and 
$
\ell_j^{[k+1]}=\Theta(\ell_j^{[k]})$.  
Consequently, these forms have rank $\nu$.

For every $j$ and $k$, one has
\[
\ell_j^{[k]}(z,\g(z))
=
(T\partial)^k\bigl(\ell_j(z,\g(z))\bigr)=(T\partial)^k\bigl(0)
=0.
\]
Hence $\ell_j^{[k]}\in\mathfrak I_\delta$, and therefore
$\ell_j^{[k]}(1,\Y)$ belongs to
${\rm ev}_1(\mathfrak I_\delta)$. Since
$$
\ell_1(1,\Y),\ldots,\ell_{\nu-\mu}(1,\Y)
$$
form a basis of ${\rm ev}_1(\mathfrak I_\delta)$ over $\mathbb K$, it
follows that every $\ell_j^{[k]}(1,\Y)$ is a linear combination of
these forms. Consequently,
\[
R_0(1,\Y),\ldots,R_{\mu+\chi_3-1}(1,\Y),
\ell_1(1,\Y),\ldots,\ell_{\nu-\mu}(1,\Y)
\]
have rank $\nu$. In view of~\eqref{eq:rho34}, the conclusion follows 
by taking
$
c_1=\chi_2$ and 
$c_2=\chi_3$.
\end{proof}

For every $n\geq n_1$, let
$L_1,\ldots,L_\mu$ be the forms provided by
Lemma~\ref{lem:formes_Pade}. We define the matrix
\[
M_n:=(a_{n;i,j})_{1\leq i,j\leq\nu}
\]
by setting
\[
a_{n;i,j}
:=
\begin{cases}
\text{the coefficient of $Y_j$ in $L_i(\Y)$},
 & 1\leq i\leq\mu,\\
\text{the coefficient of $Y_j$ in $\ell_{i-\mu}(1,\Y)$},
 & \mu<i\leq\nu.
\end{cases}
\]
By Lemma~\ref{lem:formes_Pade}, the matrix $M_n$ is invertible. Its
first $\mu$ rows correspond to forms arising from the Pad\'e
approximants, whereas its remaining $\nu-\mu$ rows are the
specializations at $z=1$ of functional relations. The sequence
$(M_n)_{n\geq n_1}$ is our sequence of preparation matrices.

\subsection{The Siegel--Shidlovskii theorem}
\label{sec:SSthm}

We are now ready to prove the Siegel--Shidlovskii theorem.

\begin{thm}[Siegel--Shidlovskii, 1959]
\label{thm:SS}
Let $f_1,\ldots,f_m$ be $E$-functions satisfying a linear differential
system, and let $\alpha\in\Q^\times$ be a regular point with respect to this system.
Then
\[
{\rm tr.deg}_{\Q}
\bigl(f_1(\alpha),\ldots,f_m(\alpha)\bigr)
=
{\rm tr.deg}_{\Q(z)}
\bigl(f_1(z),\ldots,f_m(z)\bigr).
\]
\end{thm}

\begin{proof}
We continue with the notation of
Sections~\ref{sec:notation_method} and \ref{sec:prepa_matrix}. In particular, we assume that
$\alpha=1$. By adjoining the constant function $1$ if necessary, we
may also assume that $f_m=1$. Set
\[
t:={\rm tr.deg}_{\mathbb K(z)}
\bigl(f_1(z),\ldots,f_m(z)\bigr)
\quad\mbox{ and }\quad
t_0:={\rm tr.deg}_{\mathbb K}
\bigl(f_1(1),\ldots,f_m(1)\bigr).
\]
Since $\mathbb K/\QQ$ is algebraic, these are also the transcendence
degrees appearing in the statement of the theorem. Assume, for a
contradiction, that $t_0<t$.

Since $f_m=1$, homogenization with respect to $X_m$ identifies
$\mathfrak I(1)_\delta$ with the space of polynomials in
$X_1,\ldots,X_{m-1}$ of degree at most $\delta$ which vanish at
$$
f_1(1),\ldots,f_{m-1}(1),$$ 
and identifies $(\mathfrak I_{\mathbb K(z)})_\delta$ with the space of
polynomials in $X_1,\ldots,X_{m-1}$ of degree at most $\delta$ which
vanish at
$$
f_1(z),\ldots,f_{m-1}(z).
$$ 
It follows from the Hilbert--Serre theorem (see for instance
\cite[p.~231]{FN}) that, for all sufficiently large $\delta$, the
quantities
\[
\kappa
:=
\nu-\dim_{\mathbb K}\mathfrak I(1)_\delta
\quad\mbox{ and }\quad
\mu
=
\nu-\dim_{\mathbb K(z)}
(\mathfrak I_{\mathbb K(z)})_\delta
\]
are given by polynomials in $\delta$ of degrees $t_0$ and $t$,
respectively. Since $t_0<t$, we may choose $\delta$ sufficiently large so
that
\begin{equation}
\label{eq:choice_mu}
\mu>2[\mathbb K:\QQ]\kappa.
\end{equation}

By construction,
$
{\rm ev}_1(\mathfrak I_\delta)
\subseteq
\mathfrak I(1)_\delta$. 
Moreover,
$
\ell_1(1,\Y),\ldots,\ell_{\nu-\mu}(1,\Y)$ 
form a basis of ${\rm ev}_1(\mathfrak I_\delta)$ over $\mathbb K$.
Since
$
\dim_{\mathbb K}\mathfrak I(1)_\delta=\nu-\kappa$, 
we may choose linear forms
\[
N_1(\Y),\ldots,N_\eta(\Y)
\in\mathcal O_{\mathbb K}[\Y],
\qquad
\eta:=\mu-\kappa,
\]
such that
\[
N_1(\Y),\ldots,N_\eta(\Y),
\ell_1(1,\Y),\ldots,\ell_{\nu-\mu}(1,\Y)
\]
form a basis of $\mathfrak I(1)_\delta$ over $\mathbb K$.

On the other hand,
\[
L_1(\Y),\ldots,L_\mu(\Y),
\ell_1(1,\Y),\ldots,\ell_{\nu-\mu}(1,\Y)
\]
form a basis of the $\mathbb K$-vector space of all linear forms in
$\Y$. The images of $L_1,\ldots,L_\mu$ therefore span the quotient of
this space by $\mathfrak I(1)_\delta$, which has dimension $\kappa$.
After reindexing the forms $L_i$, we may consequently assume that
\begin{equation}
\label{eq:forms}
L_1(\Y),\ldots,L_\kappa(\Y),
N_1(\Y),\ldots,N_\eta(\Y),
\ell_1(1,\Y),\ldots,\ell_{\nu-\mu}(1,\Y)
\end{equation}
form a basis of the space of all linear forms in $\Y$.

Let
$
M_n'\in\mathcal M_\nu(\mathcal O_{\mathbb K})$ 
be the matrix whose rows are the coefficient vectors, with respect to
$Y_1,\ldots,Y_\nu$, of the forms in~\eqref{eq:forms}. This will be our forbidden matrix. 

Since
\[N_i(\bw)=0
\qquad
(1\leq i\leq\eta),
\quad\mbox{ and }\quad 
\ell_i(1,\bw)=0
\qquad
(1\leq i\leq\nu-\mu),
\]
we have
\begin{equation}\label{SSForb}
M_n'\bw
=
{}^{\rm t}
\bigl(
L_1(\bw),\ldots,L_\kappa(\bw),0,\ldots,0
\bigr)
=:\boldsymbol\varepsilon_n.
\end{equation}
By Lemma~\ref{lem:formes_Pade},
$
\| \boldsymbol\varepsilon_n\|_\infty
\leq
\rho_1n^{-(\mu-1-\varepsilon)n}$. 

Let $C_n$ denote the maximum of the absolute values of the entries in
the first $\kappa$ columns of $(M_n')^{-1}$. Since only the first
$\kappa$ coordinates of $\boldsymbol\varepsilon_n$ may be non-zero, we obtain
\begin{equation}
\label{eq:borne_w}
\|\bw\|_\infty
=
\|(M_n')^{-1}\boldsymbol\varepsilon_n\|_\infty
\leq
\kappa C_n\| \boldsymbol\varepsilon_n\|_\infty
\leq
\kappa\rho_1C_n
n^{-(\mu-1-\varepsilon)n}.
\end{equation}

We now estimate $C_n$ using the comatrix formula. Set
$
H_n:=\rho_2n^{(1+\varepsilon)n}$. 
Then $H_n$ is an upper bound for the houses of
$L_1,\ldots,L_\kappa$. Let $H$ be an upper bound for the houses of the
coefficients of
\[
N_1,\ldots,N_\eta,
\ell_1(1,\Y),\ldots,\ell_{\nu-\mu}(1,\Y).
\]
These forms are independent of $n$, and hence so is $H$.

Since
$
\det M_n'\in\mathcal O_{\mathbb K}\setminus\{0\}$, 
Liouville's inequality gives
\[
|\det M_n'|
\geq
\house{\det M_n'}^{\,1-[\mathbb K:\QQ]}.
\]
Moreover,
\[
\house{\det M_n'}
\leq
\nu!H_n^\kappa H^{\nu-\kappa}.
\]
Consequently, there exists a positive constant $\rho_3$, independent
of $n$, such that
\begin{equation}
\label{eq:liouville_det}
|\det M_n'|
\geq
\rho_3H_n^{-\kappa([\mathbb K:\QQ]-1)}.
\end{equation}

An entry in one of the first $\kappa$ columns of
$\operatorname{Com}(M_n')^{\rm t}$ is a cofactor obtained by deleting
one of the first $\kappa$ rows of $M_n'$. It therefore involves at
most $\kappa-1$ rows whose coefficients have house bounded by $H_n$.
Hence its absolute value is at most
$
(\nu-1)!H_n^{\kappa-1}H^{\nu-\kappa}$. 
Combining this estimate with~\eqref{eq:liouville_det}, we obtain a
positive constant $\rho_4$, independent of $n$, such that
\begin{equation}
\label{eq:bound_Cn}
C_n
\leq
\rho_4H_n^{\kappa[\mathbb K:\QQ]-1}.
\end{equation}

Combining~\eqref{eq:borne_w} and~\eqref{eq:bound_Cn}, and absorbing
all constants independent of $n$ into a positive constant $\rho_5$,
we obtain
\begin{align}
\|(M_n')^{-1}\boldsymbol\varepsilon_n\|_\infty 
&\leq
\rho_5
n^{(1+\varepsilon)n(\kappa[\mathbb K:\QQ]-1)}
n^{-(\mu-1-\varepsilon)n}
\nonumber\\
&=
\rho_5
n^{((1+\varepsilon)\kappa[\mathbb K:\QQ]-\mu)n}.
\label{eq:contradiction}
\end{align}

By~\eqref{eq:choice_mu} and the assumption $0<\varepsilon<1$, we have
$
(1+\varepsilon)\kappa[\mathbb K:\QQ]-\mu<0$. 
The right-hand side of~\eqref{eq:contradiction} therefore tends to
zero as $n$ tends to infinity. Since $\bw\neq0$, this contradicts
\eqref{eq:borne_w}.  
Thus $t_0\geq t$. The reverse inequality is immediate by
specialization of algebraic relations. We conclude that $t_0=t$,
which completes the proof.
\end{proof}

\subsection{The lifting theorem}\label{sec:Beukers}

We now prove Beukers' lifting theorem.

\begin{thm}[Beukers, 2006]\label{thm:Beukers}
Under the assumptions of the Siegel--Shidlovskii Theorem, every
homogeneous polynomial $P\in\Q[X_1,\ldots,X_m]$ such that
\[
P(f_1(\alpha),\ldots,f_m(\alpha))=0
\]
can be lifted to a polynomial
$Q\in\Q[z,X_1,\ldots,X_m]$, homogeneous in
$X_1,\ldots,X_m$, such that
\[
Q(z,f_1(z),\ldots,f_m(z))=0
\quad\mbox{and}\quad
Q(\alpha,X_1,\ldots,X_m)
=
P(X_1,\ldots,X_m).
\]
\end{thm}

\begin{proof} We continue with the notation of
	Sections~\ref{sec:notation_method} and \ref{sec:prepa_matrix}. In particular, we assume that
	$\alpha=1$.
  It suffices to prove that
\begin{equation}\label{eq:egalite_ideals}
\mathfrak I(1)
=
{\rm ev}_1(\mathfrak I).
\end{equation} 
By definition, $
{\rm ev}_1(\mathfrak I)
\subset
\mathfrak I(1)$. Assume, for a contradiction, that
$
{\rm ev}_1(\mathfrak I)
\subsetneq
\mathfrak I(1)$. 

By~{\rm (A2)}, reproved in Appendix~\ref{AppendiceA}, the ideal
${\rm ev}_1(\mathfrak I)$ is prime. The ideal $\mathfrak I(1)$ is also
prime, since it is the homogeneous ideal of algebraic relations among
$f_1(1),\ldots,f_m(1)$. 
Therefore, if $h$ and $h_0$ denote their
respective heights (see \cite[p.\ 193]{ZS} for a definition), then $h_0>h$. Recall that the basis of $\mathfrak I_\delta$ was chosen so that its
specialization at $z=1$ remains linearly independent. In particular,
$
\dim_{\mathbb K}{\rm ev}_1(\mathfrak I)_\delta
=
\dim_{\mathbb K(z)}\mathfrak I_\delta$.

By the Hilbert--Serre Theorem, for sufficiently large $\delta$, the numbers
\[
\mu
=
\nu-\dim_{\mathbb K(z)}\mathfrak I_\delta=\nu-\dim_{\mathbb K}{\rm ev}_1(\mathfrak I)_\delta,
\quad \text{and} \quad
\kappa
:=
\nu-\dim_{\mathbb K}\mathfrak I(1)_\delta.
\]
agree with polynomials of respective degrees
$m-h-1$ and $m-h_0-1$. Since $h_0>h$, we may choose $\delta$
sufficiently large so that
$
\mu>2[\mathbb K:\mathbb Q]\kappa$.

We now argue as in the proof of Theorem~\ref{thm:SS}. Choose a
complement of ${\rm ev}_1(\mathfrak I)_\delta$ in
$\mathfrak I(1)_\delta$. Its dimension is $\eta:=\mu-\kappa$, and a basis of
this complement provides $\eta$ additional putative forms. These
forms may be substituted for $\eta$ of the Pad\'e forms occurring
in the fixed preparation matrix, leaving only $\kappa$ Pad\'e forms. 
Exactly the same comatrix argument as in the proof of
Theorem~\ref{thm:SS}, together with
$
\mu>2[\mathbb K:\mathbb Q]\kappa$, 
shows that the resulting matrix would be a forbidden matrix. This
contradiction proves~\eqref{eq:egalite_ideals}, and the theorem follows.
\end{proof}
\subsection{Algebraic independence measures}
\label{sec:measures}

We now prove Theorem~\ref{thm:measure_system} stated in the introduction.

We continue with the notation of Sections~\ref{sec:notation_method} and \ref{sec:prepa_matrix} and, in particular, we assume $\alpha=1$.
Assume for a contradiction that the conclusion of Theorem~\ref{thm:measure_system} fails. Then, for arbitrarily large real numbers $\gamma$ there exist an integer $d$ and infinitely many polynomials $P\in\mathbb Z[X_1,\ldots,X_m]$ with 
degree at most $d$ and height at most $H$ such that
\begin{equation}\label{eq:contradiction_mesure}
0 <  \vert P(f_1(1),\ldots,f_m(1))\vert < H^{-\gamma d^t},
\end{equation}
where
$
t:={\rm tr.deg}_{\mathbb K(z)}
\bigl(f_1(z),\ldots,f_m(z)\bigr)$.
Fix a parameter $\gamma$ large enough, an integer $d$ associated, and a polynomial $P$ satisfying \eqref{eq:contradiction_mesure}, with $H$ large enough, to be specified later. By adjoining the
constant function $1$ and homogenizing $P$, we may assume that $f_m=1$
and that $P$ is homogeneous. The fact that $P(f_1(1),\ldots,f_m(1))\neq 0$ yields $P\notin\mathfrak I(1)$. Let $\mathfrak P$ denote the homogeneous
ideal of $\mathbb K[\X]$ generated by $P$ and $\mathfrak I(1)$.

\subsubsection{Preliminary dimension estimates}
\label{sec:dimension}

Let $\tau$ be a positive integer, to be chosen later,
and set
$
\delta=d\tau$. 
Let $\nu$ and $\mu$ be the integers associated with $\delta$ in
Section~\ref{sec:prepa_matrix}. We set
\[
\kappa
:=
\nu-\dim_{\mathbb K}\mathfrak P_\delta
\quad\mbox{ and }\quad
\eta
:=
\dim_{\mathbb K}\mathfrak P_\delta
-
\dim_{\mathbb K(z)}
(\mathfrak I_{\mathbb K(z)})_\delta
=
\mu-\kappa.
\]

\begin{lem}
\label{lem:7-8}
There exist a positive integer $\tau_1$ and a positive real number
$\lambda$, both independent of $\gamma,d$ and $P$, such that, for every
$\tau\geq\tau_1$,
\[
\mu > 2[\mathbb K:\mathbb Q]\kappa
\quad\mbox{ and }\quad 
\mu\leq\lambda\tau^td^t.
\] 
\end{lem}

\begin{proof}
The second estimate follows from the fact that, for all sufficiently
large $\delta$, the integer $\mu$ is given by a polynomial in
$\delta=d\tau$ of degree $t$, as observed in the proof of
Theorem~\ref{thm:SS}.

Set
$
\mathcal A
:=
\mathbb K[\X]/\mathfrak I(1)$, 
and let $\mathcal A_\delta$ denote its homogeneous component of degree
$\delta$. Since $\mathfrak I(1)$ is prime, $\mathcal A$ is an integral
domain. Let $\varphi$ be its Hilbert polynomial. For all sufficiently
large $\tau$, uniformly with respect to $d\geq1$, we have
\[
\dim_{\mathbb K}\mathcal A_{\tau d}
=
\varphi(\tau d)
\quad\mbox{ and }\quad 
\dim_{\mathbb K}\mathcal A_{(\tau-1)d}
=
\varphi((\tau-1)d).
\]

Let $\overline P$ denote the image of $P$ in $\mathcal A$. Since
$\overline P\neq0$ and $\mathcal A$ is an integral domain,
multiplication by $\overline P$ induces an isomorphism
$
\mathcal A_{(\tau-1)d}
\longrightarrow
(\overline P\mathcal A)_{\tau d}$. 
Consequently,
$
\dim_{\mathbb K}(\overline P\mathcal A)_{\tau d}
=
\varphi((\tau-1)d)$. 
Since
\[
\frac{\varphi((\tau-1)d)}{\varphi(\tau d)}
\longrightarrow1
\]
as $\tau$ tends to infinity, uniformly with respect to $d\geq1$, there
exists $\tau_1$, independent of $P$, such that
\[
\frac{
\dim_{\mathbb K}(\overline P\mathcal A)_\delta
}{
\dim_{\mathbb K}\mathcal A_\delta
}
>
1-\frac{1}{2[\mathbb K:\QQ]}
\]
whenever $\tau\geq\tau_1$. 
Now
$
\kappa
=
\dim_{\mathbb K}\mathcal A_\delta
-
\dim_{\mathbb K}(\overline P\mathcal A)_\delta$. 
Moreover, since
$
{\rm ev}_1(\mathfrak I_\delta)
\subseteq
\mathfrak I(1)_\delta$, 
we have
$
\dim_{\mathbb K}\mathcal A_\delta\leq\mu$. 
It follows that
$
2[\mathbb K:\QQ]\kappa
<
\mu$, 
as desired.
\end{proof}

We henceforth fix $\tau=\tau_1$ and set
$
\varepsilon=\frac12$. 
Writing $D:=[\mathbb K:\QQ]$, Lemma~\ref{lem:7-8} gives
\begin{equation}
\label{eq:gap_measure}
\mu-(1+\varepsilon)D\kappa>\frac{\mu}{4}.
\end{equation}
We also assume, without loss of generality, that $\gamma$ is chosen large enough so that \begin{equation}\label{eq:choice_gamma}
\gamma \geq 16D\lambda \tau^t. 
\end{equation}

\subsubsection{The putative forms}

Let $(M_n)_{n\geq n_1}$ be the sequence of preparation matrices
constructed in Section~\ref{sec:prepa_matrix} for $\delta=d\tau$ and
$\varepsilon=1/2$. Recall that its rows are the coefficient vectors of
\[
L_1(\Y),\ldots,L_\mu(\Y),
\ell_1(1,\Y),\ldots,\ell_{\nu-\mu}(1,\Y).
\]
We set
\[
\bw
:=
{}^{\rm t}(g_1(1),\ldots,g_\nu(1)),
\]
where $g_1,\ldots,g_\nu$ are the monomials of degree $\delta$ in
$f_1,\ldots,f_m$. 
The quotient
$
\mathfrak P_\delta/\mathfrak I(1)_\delta$ 
identifies with  the degree-$\delta$ component of the principal ideal generated by
the image of $P$ in $\mathcal A$. We may therefore choose monomials
$A_1(\X),\ldots,A_\eta(\X)$ of degree $\delta-d$ such that the forms
\[
N_i(\Y):=A_i(\X)P(\X),
\qquad
1\leq i\leq\eta,
\]
together with 
$
\ell_1(1,\Y),\ldots,\ell_{\nu-\mu}(1,\Y)$ are linearly independent\footnote{Note that it follows from Beukers's lifting theorem that these forms span $\mathfrak P_\delta$, 
	an observation that we shall not need here.}, since the latter belong to 
$\mathfrak I(1)_\delta$. 
These are our putative forms. In particular,
\[
\house{N_i}\leq H,
\qquad
1\leq i\leq\eta.
\]

\subsubsection{End of the proof of Theorem~\ref{thm:measure_system}}
After reindexing the Pad\'e forms, we may assume that
\begin{equation}
\label{eq:forms2}
L_1(\Y),\ldots,L_\kappa(\Y),
N_1(\Y),\ldots,N_\eta(\Y),
\ell_1(1,\Y),\ldots,\ell_{\nu-\mu}(1,\Y)
\end{equation}
form a basis of the space of linear forms in $\Y$. 
Let
$
M_n'\in\mathcal M_\nu(\mathcal O_{\mathbb K})$ 
be the matrix whose rows are the coefficient vectors of the forms
in~\eqref{eq:forms2}. Then
\[
M_n'\bw
=
{}^{\rm t}\bigl(\e_n,\e_n',\boldsymbol 0\bigr)=:\boldsymbol\varepsilon_n,
\]
where
$
\e_n
:=
\bigl(L_1(\bw),\ldots,L_\kappa(\bw)\bigr)
$
and
$
\e_n'
:=
\bigl(N_1(\bw),\ldots,N_\eta(\bw)\bigr)$. 
By Lemma~\ref{lem:formes_Pade},
$
\|\e_n\|_\infty
\leq
\rho_1n^{-(\mu-1-\varepsilon)n}$. 
Moreover, since $N_i=A_iP$, we infer from \eqref{eq:contradiction_mesure} that there exists 
a positive real number
$\theta$, depending only on $f_1(1),\ldots,f_m(1)$, such that
\begin{equation}
\label{eq:mesure_Ni}
\|\e_n'\|_\infty
\leq
\theta^{\delta-d}
 H^{-\gamma d^t}.
\end{equation}

Let $C_n$ be the maximum of the absolute values of the entries in the
first $\kappa$ columns of $(M_n')^{-1}$, and let $C_n'$ be the
corresponding maximum over the following $\eta$ columns. Up to factors
depending only on $\delta$, which are absorbed into the constants, we
obtain
\begin{equation}
\label{eq:borne_w2}
\|(M_n')^{-1}\boldsymbol \varepsilon_n\|_\infty
\leq
\rho_1C_n n^{-(\mu-1-\varepsilon)n}
+
C_n' H^{-\gamma d^t}.
\end{equation}

Set
$
H_n:=\rho_2n^{(1+\varepsilon)n}$. 
The first $\kappa$ rows of $M_n'$ have house at most $H_n$, the
following $\eta$ rows have house at most $H$, and the remaining
rows have house bounded independently of $n$ and $P$. Since
$\det M_n'$ is a non-zero algebraic integer, Liouville's inequality and
the comatrix formula give
\[
|\det M_n'|^{-1}
\leq
\rho_3
H_n^{\kappa(D-1)}
H^{\eta(D-1)}
\]
for some positive constant $\rho_3$ independent of $n$ and $H$. 
The cofactors corresponding to the first $\kappa$ columns of
$(M_n')^{-1}$ involve at most $\kappa-1$ rows of house $H_n$ and $\eta$
rows of house $H$, whereas those corresponding to the following
$\eta$ columns involve at most $\kappa$ rows of house $H_n$ and
$\eta-1$ rows of house $H$. Hence
\begin{align}
C_n
&\leq
\rho_4
n^{(1+\varepsilon)n(\kappa D-1)}
H^{\eta D},
\label{eq:bound_Cn2}\\
C_n'
&\leq
\rho_4
n^{(1+\varepsilon)n\kappa D}
H^{\eta D-1},
\nonumber
\end{align}
where $\rho_4$ is independent of $n$ and $H$. Combining
\eqref{eq:borne_w2} and~\eqref{eq:bound_Cn2}, we obtain
\[
\|(M_n')^{-1}\boldsymbol \varepsilon_n\|_\infty
\leq
\rho_1\rho_4
n^{-n(\mu-(1+\varepsilon)\kappa D)}
H^{\eta D}
+
\rho_4
n^{(1+\varepsilon)n\kappa D}
H^{\eta D-\gamma d^t-1}.
\]
Let $n$ be the unique integer such that
$$
\rho_1^{-1}n^{n\mu} \leq H^{\gamma d^t+1} < \rho_1^{-1} (n+1)^{(n+1)\mu}\,.
$$
We may assume $H$ is large enough to ensure that $n \geq n_1$. Therefore we have
\[
\|(M_n')^{-1}\boldsymbol \varepsilon_n\|_\infty
\leq
2\rho_1\rho_4
n^{-n(\mu-(1+\varepsilon)\kappa D)}
H^{\eta D}. 
\]
Our choice of $n$, Lemma~\ref{lem:7-8}, \eqref{eq:choice_gamma}, and the rough estimates $\gamma d^t + 1\leq 2\gamma d^t$, and $(n+1)^{(n+1)}\leq n^{2n}$, imply that
$$
H \leq \rho_5 n^{4n\mu  \gamma^{-1}d^{-t}} \leq \rho_5n^{4n\lambda \tau^t \gamma^{-1}}\leq \rho_5 n^{n/(4D)},$$
for some positive real number $\rho_5$.
Therefore
\begin{equation}\label{eq:majo_inversematrix}
\|(M_n')^{-1}\boldsymbol \varepsilon_n\|_\infty
\leq \rho_6
n^{-n(\mu-(1+\varepsilon)\kappa D-\eta/4)}\,.
\end{equation}
Using \eqref{eq:gap_measure} and the fact that $\eta \leq \mu$, we have $\mu-(1+\varepsilon)\kappa D-\eta/4>0$. 
Since $n$ tends to infinity with $H$, if $H$ is large enough, the quantity on the right-hand side of \eqref{eq:majo_inversematrix} is smaller than $\|\bw\|_\infty$, which contradicts that $(M_n')^{-1}\boldsymbol \varepsilon_n=\bw$. Thus, $M_n'$ is a forbidden matrix, which ends the proof.
\qed

\subsubsection{Remark on effectivity}\label{sec:effectivity}

Since Lang's work~\cite{La62}, several authors have investigated the
effectivity of the constants appearing in algebraic independence
measures for values of $E$-functions. When
$f_1,\ldots,f_m$ are algebraically independent, it follows from
\cite[Thm.~1]{Gor90} that one may take
\begin{equation}\label{eq:bounds_ci}
C_1
=
\exp\left(
-\exp\left(
C_3d^{2m}\log(d+1)
\right)
\right)
\ \text{and}\ 
C_2
=
\sqrt{m}\,4^m[\mathbb K:\mathbb Q]^{m+1},
\end{equation}
where $C_3>0$ depends only on $f_1,\ldots,f_m$ and $\alpha$.
Obtaining an effective value for $C_3$ remained a long-standing
problem, since it depends on the Shidlovskii constant arising in
Shidlovskii's lemma; see, for instance,
\cite[Lem.~5.2, p.~221]{FN}. This problem was eventually solved through
the work of Bertrand and Beukers~\cite{BB85} and of Bertrand, Chirskii,
and Yebbou~\cite{BCY04}.

The constants occurring in the more general setting of
Theorem~\ref{thm:measure_system} can likewise be bounded effectively.
We omit the details, since making the constants explicit is standard, 
although rather cumbersome. Let us only indicate the additional
ingredients that arise in the present situation.

The constants $\lambda$ and
$\max_i\house{\ell_i(1,\Y)}$ do not appear in the classical argument,
where $f_1,\ldots,f_m$ are assumed to be algebraically independent.
Moreover, in that case, the parameter $\tau_1$ depends only on $m$ and
$[\mathbb K:\mathbb Q]$, and this dependence can be made explicit. In
our more general setting, these quantities also depend on the ideal of
functional relations $\mathfrak I$. The main difficulty is therefore to bound the degrees and the heights of a generating set for
$\mathfrak I$; see the discussion in~\cite{AF24}.

Once these constants and the parameter $n_1$  arising from Shidlovskii's lemma have been bounded the remaining constants can be estimated by means of the results of~\cite{Sh79}.

We finally record, without proof, that the constants in
Theorem~\ref{thm:measure_system} may be chosen in the form
\[
C_1
=
\exp\left(
-\exp\left(
C_3[\mathbb K:\mathbb Q]^{2t}
d^{2t}
\log\bigl([\mathbb K:\mathbb Q]d+1\bigr)
\right)
\right)
\quad\text{and}\quad
C_2
=
C_4[\mathbb K:\mathbb Q]^{t+1},
\]
where $C_3,C_4>0$ depend only on
$f_1,\ldots,f_m$ and $\alpha$, but not on $d$,
$[\mathbb K:\mathbb Q]$, or $t$.

\section{The Siegel--Shidlovskii method at singular points}\label{sec:Andre}

Theorems~\ref{thm:SS}, \ref{thm:Beukers}, and
\ref{thm:measure_system} concern values of $E$-functions at a point
that is regular for a differential system satisfied by these
functions. This regularity assumption allows the derivatives of the
initial Pad\'e approximant to be rewritten as linear combinations of
the original functions, thereby ensuring the existence of a sequence
of preparation matrices.

At a singular point, such a rewriting is no longer available in
general. Nevertheless, the successive derivatives of the Pad\'e
approximant can still be interpreted as linear forms in a larger
family of functions. This observation extends the
Siegel--Shidlovskii method to points that may a priori be singular.
We do not pursue this framework in full generality here, but show that
it is already sufficient to give a new proof of Andr\'e's theorem~{\rm (A1)}, which we restate below in the following slightly different form.

\begin{thm}\label{thm:Andre}
Every $E$-function satisfies a linear differential equation whose 
singularities belong to $\{0,\infty\}$.
\end{thm}

We
then combine Theorem~\ref{thm:Andre} with
Theorem~\ref{thm:measure_system} to prove
Corollary~\ref{thm:measure_sans_system}. We conclude this section by discussing some
further consequences of Theorem~\ref{thm:Andre}. In
Appendix~\ref{AppendiceB}, we will use Theorem~\ref{thm:Andre} to obtain a refinement of {\rm (A1)}, as well as a consequence concerning the
minimal differential equation of a $G$-function.

\subsection{The rational case}

We begin by proving a result for $E$-functions with
rational coefficients, from which the general case will follow.

\begin{thm}\label{thm:Andre_surQ}
Let $f\in\mathbb Q[[z]]$ be an $E$-function. Then $f$ satisfies a
linear differential equation having $1$ as an ordinary point.
\end{thm}

\subsubsection{Preliminary setting}

If $f=0$, there is nothing to prove, so we assume that $f\neq0$. Let
$i_0$ be the least non-negative integer such that $f^{(i_0)}(1)\neq0$,
and let $\mu\geq1$ be the order of the minimal differential equation
satisfied by $f$. Let $\varepsilon>0$ be a sufficiently small real number satisfying
$\varepsilon\mu<1$. We consider the infinite family of $E$-functions
$
\f:=(f^{(i)}(z))_{i\geq0}\subset\mathbb Q[[z]]$ 
and the corresponding family of indeterminates
$\X:=(X_i)_{i\geq0}$.

Since we work over $\mathbb Q$, it is enough to consider linear forms
in $f$ and its derivatives. However, we do not fix in advance the
number $\nu$ of derivatives involved, and allow it to depend on $n$.
Consequently, the size of the preparation matrix $M_n$ will also
depend on $n$. For fixed $n$, the rows of $M_n$ are divided into two
blocks with, respectively, $\mu$ and $\nu-\mu$ rows:
\begin{itemize}
\item the first block consists of the specializations at $z=1$ of
certain derivatives of a Pad\'e approximant involving
$f,f',\ldots,f^{(\mu-1)}$;

\item the second block consists of the specializations at $z=1$ of
differential relations involving
$f,f',\ldots,f^{(\nu-1)}$.
\end{itemize}

\subsubsection{Forms arising from the differential equations}

For every $i\geq\mu$, we let $d_i$ denote the minimum of the degree of $a_i(z)$ among all differential relations 
\[
a_{0}(z)f(z)+a_{1}(z)f'(z)+\cdots
+a_{i}(z)f^{(i)}(z)=0,
\qquad a_{0},\ldots,a_i\in\mathbb Z[z],\ a_i \neq 0.
\]
Since differentiating a relation of order
$i$ produces a relation of order $i+1$ with the same leading
coefficient, the sequence $(d_i)_{i\geq\mu}$ is non-increasing. Let
$i_1$ be such that this sequence is constant for $i\geq i_1$. 

For each $i \geq \mu$,
we choose a differential relation
\[
a_{i,0}(z)f(z)+a_{i,1}(z)f'(z)+\cdots
+a_{i,i}(z)f^{(i)}(z)=0,
\qquad a_{i,j}\in\mathbb Z[z], 
\]
with  $\deg a_{i,i}=d_i$ by proceeding as follows. When $i \leq i_1$, we select an arbitrary differential relation for which $\deg a_{i,i}=d_i$. When $i > i_1$, we choose the differential relation obtained by differentiating $i-i_1$ times the relation of order $i_1$.

Set
\[
\ell_i(z,\X)
:=
a_{i,0}(z)X_0+\cdots+a_{i,i}(z)X_i.
\]
Thus $\ell_i(z,\f)=0$, and, for any
$i>i_1$,
\[
\ell_{i}=\partial^{i-i_1}\ell_{i_1}
\]
where $\partial$ is the derivation of $\mathbb Z[z,\X]$ defined by
\begin{equation}\label{eq:derivation}
\partial X_i:=X_{i+1}
\qquad\text{and}\qquad
\partial a(z):=a'(z).
\end{equation}

We claim that the infinite family of linear forms
$\ell_i(1,\X)$, $i\geq\mu$, is linearly independent. Suppose, to the
contrary, that there exist integers
$\lambda_1,\ldots,\lambda_k$, with $\lambda_k\neq0$, and indices
$
\mu\leq j_1<\cdots<j_k$ 
such that
$
\lambda_1\ell_{j_1}(1,\X)+\cdots
+\lambda_k\ell_{j_k}(1,\X)=0$. 
The coefficients of
\[
\lambda_1\ell_{j_1}(z,\X)+\cdots
+\lambda_k\ell_{j_k}(z,\X)
\]
are then divisible by $z-1$. Dividing by $z-1$ yields a linear 
differential relation of order $j_k$ with polynomial coefficients satisfied by $f$, whose leading
coefficient is
$
\frac{\lambda_k a_{j_k,j_k}(z)}{z-1}$. This contradicts
the minimality of the degree of $a_{j_k,j_k}$.

The houses of these linear forms may grow with $i$, but at most
polynomially. 

\begin{lem}\label{lem:rel_func}
There exist a positive constant $c_0$ and a non-negative integer
$\delta$, independent of $i$, such that
\[
\house{\ell_i(1,\X)}\leq c_0i^\delta
\]
for every $i\geq\mu$.
\end{lem}

\begin{proof}
	Enlarging $c_0$ if necessary, we only need to establish the bound when $i\geq i_1$. Fix an integer $i \geq i_1$.
Write
\[
\ell_{i_1}(z,\X)
:=
a_0(z)X_0+\cdots+a_{i_1}(z)X_{i_1},
\qquad a_0,\ldots,a_{i_1}\in\mathbb Z[z],
\]
and set
$
\delta:=\max_{0\leq s\leq i_1}\deg a_s$. 
By the Leibniz rule,
\[
\partial^{i-i_1}(a_sX_s)
=
\sum_{j=0}^{\delta}
\binom{i-i_1}{j}
a_s^{(j)}(z)X_{s+i-i_1-j},
\]
where the terms with $j>i-i_1$ are understood to be zero. 
For convenience, we extend the family
$a_0,\ldots,a_{i_1}$ by setting
$
a_s:=0$  for every $s\notin\{0,\ldots,i_1\}$. 
Then the coefficient of $X_{i-k}$ in $\ell_i=\partial^{i-i_1}\ell_{i_1}$ is
\begin{equation}\label{eq:Leibniz}
\sum_{j=0}^{\delta}
\binom{i-i_1}{j}
a_{i_1+j-k}^{(j)}(z).
\end{equation}

Set
\[
H_0
:=
\max{\substack{0\leq s\leq i_1\\0\leq j\leq\delta}}
\left|a_s^{(j)}(1)\right|.
\]
The absolute value of every coefficient of
$\ell_i(1,\X)$ is at most
\[
H_0\sum_{j=0}^{\delta}\binom{i-i_1}{j}
\leq c_0i^\delta
\]
for a suitable positive constant $c_0$.
\end{proof}

As in the proof of
Lemma~\ref{lem:rel_func}, we write
\[
\ell_{i_1}(z,\X)
:=
a_0(z)X_0+\cdots+a_{i_1}(z)X_{i_1}.
\]
Set
\begin{equation}\label{eq:def_eta}
\eta
:=
\min_{0\leq k\leq i_1}
\left\{
k+\operatorname{ord}_{z=1}a_{i_1-k}(z)
\right\}.
\end{equation}
The polynomials $a_0,\ldots,a_{i_1}$ are relatively prime, since
otherwise their common divisor could be removed, contradicting the
minimality of $\deg a_{i_1,i_1}$. In particular, at least one of them
does not vanish at $1$, and hence $\eta\leq i_1$.

Recall the convention that $a_s=0$ for
$s\notin\{0,\ldots,i_1\}$. For every $i\geq i_1$,
Formula~\eqref{eq:Leibniz} gives
\[
\ell_i(z,\X)
=
\sum_{k=0}^{i}
\sum_{j=0}^{k}
\binom{i-i_1}{j}
a_{i_1+j-k}^{(j)}(z)X_{i-k}.
\]
By the definition of $\eta$, we have
$
a_{i_1+j-k}^{(j)}(1)=0$ 
for every $k<\eta$ and every $j\leq k$. Therefore
$\ell_i(1,\X)$ is a linear form in
$
X_0,\ldots,X_{i-\eta}$. 
The coefficient of $X_{i-\eta}$ is
\begin{equation}\label{eq:polynome_binomial}
\sum_{j=0}^{\eta}
\binom{i-i_1}{j}
a_{i_1+j-\eta}^{(j)}(1).
\end{equation}

Choose $k$ realizing the minimum in~\eqref{eq:def_eta}, and set
$
j:=\operatorname{ord}_{z=1}a_{i_1-k}(z)$. 
Then $\eta=k+j$ and
$
a_{i_1+j-\eta}^{(j)}(1)
=
a_{i_1-k}^{(j)}(1)\neq0$. 
Consequently, not all the coefficients
$
a_{i_1+j-\eta}^{(j)}(1)$,
$0\leq j\leq\eta$, 
vanish. Since the polynomials
$\binom{i-i_1}{j}$, $0\leq j\leq\eta$, are linearly independent,
Expression~\eqref{eq:polynome_binomial} is a non-zero polynomial in
$i$. There thus exists an integer $i_2\geq i_1$ such that the
coefficient of $X_{i-\eta}$ in $\ell_i(1,\X)$ is non-zero for every
$i\geq i_2$. Enlarging $i_2$ if necessary, we may assume that
$i_2\geq i_0$.


\subsubsection{Forms arising from Pad\'e approximants}

By the minimality of the differential equation satisfied by $f$, the
functions
$
f,f',\ldots,f^{(\mu-1)}$ 
are linearly independent over $\mathbb Q(z)$. Moreover, they form a
solution of the companion system associated with this equation.

Recall that $\varepsilon>0$ was chosen so that $\varepsilon\mu<1$. We
apply Proposition~\ref{prop:SiegelShidlovskii} to these $\mu$ functions
with
\[
\varepsilon_0:=\frac{\varepsilon}{6(\mu+1)}.
\]
Let $n\geq n_1$, where $n_1$ is given by that proposition. There exist
a positive real number $\rho$, independent of $n$, and polynomials
$P_1,\ldots,P_\mu\in\mathbb Z[z]$ of degree at most $n$ such that:
\begin{itemize}
\item[-]
$
\house{P_i}\leq \rho n^{(1+\varepsilon_0)n}
\qquad (1\leq i\leq\mu)$;

\item[-] setting
$
R(z):=P_1(z)f(z)+\cdots+P_\mu(z)f^{(\mu-1)}(z)$, 
we have
$$
\val R(z)\geq \mu(n+1)-\varepsilon_0n-1;$$

\item[-] the functions
$
R,R',\ldots,R^{(\mu-1)}$ 
are linearly independent over $\mathbb Q(z)$.
\end{itemize}

We now define a sequence $(R_k)_{k\geq0}$ of linear forms in
$\X=(X_i)_{i\geq0}$ by
\[
R_0(z,\X)
:=
P_1(z)X_0+\cdots+P_\mu(z)X_{\mu-1}
\]
and
\[
R_k(z,\X):=\partial^kR_0(z,\X),
\]
where $\partial$ denotes the derivation of $\mathbb Z[z,\X]$ defined
by \eqref{eq:derivation}. Then
$
R_k(z,\f)=R^{(k)}(z)$ 
for every $k\geq0$, and $R_k$ involves only the variables
$
X_0,\ldots,X_{\mu+k-1}$. 

\subsubsection{The preparation matrix}

We first prove the following statement.

\begin{lem}\label{lem:sing_good_forms}
There exist a positive real number $c_1$, independent of $n$ and
$\varepsilon$, two integers $\nu\geq i_2$ and $c_2\geq 0$ satisfying
$
c_2+\mu\leq \nu\leq \varepsilon_0 n+c_1$, 
positive real numbers $\rho_1,\rho_2$, and $\mu$ forms
$L_1(\X),\ldots,L_\mu(\X)$ among
$
R_0(1,\X),\ldots,R_{\mu+c_2-1}(1,\X)$ 
such that:

\medskip

\begin{itemize}
\item[\rm(a)]
$
|L_i(\f(1))|
\leq \rho_1n^{-(\mu-1-\varepsilon)n}
\qquad (1\leq i\leq\mu)$;

\medskip

\item[\rm(b)]
$
\house{L_i}\leq \rho_2n^{(1+\varepsilon)n}
\qquad (1\leq i\leq\mu)$;

\medskip

\item[\rm(c)]
the forms
$
L_1(\X),\ldots,L_\mu(\X),
\ell_\mu(1,\X),\ldots,\ell_{\nu-1}(1,\X)$ 
span the $\mathbb Q$-vector space of linear forms in
$X_0,\ldots,X_{\nu-\eta-1}$.
\end{itemize}
\end{lem}

\begin{proof}
Assertions~(a) and~(b) follow directly from the construction of the forms
$R_k(z,\X)$. To
prove~(c), we first show that there exist $\nu$ and $c_2$ such that
$X_0,\ldots,X_{\nu-\eta-1}$ belong to the $\mathbb Q$-vector space
spanned by
\[
R_0(1,\X),\ldots,R_{\mu+c_2-1}(1,\X),
\ell_\mu(1,\X),\ldots,\ell_{\nu-1}(1,\X).
\]

Since $\nu$ depends on $n$, the argument differs slightly from that
used in the proof of Lemma~\ref{lem:formes_Pade}. Set
$
i_3:=\max\{i_2+1,2\mu\}$, 
and let
\[
\Gamma(z):=(\gamma_{i,j}(z))_{0\leq i,j\leq i_3-1}
\]
be defined by
\[
\gamma_{i,j}(z):=
\begin{cases}
\text{the coefficient of $X_j$ in $R_i(z,\X)$},
& 0\leq i<\mu,\\
\text{the coefficient of $X_j$ in $\ell_i(z,\X)$},
& \mu\leq i<i_3.
\end{cases}
\]
Arguing as in the proof of Lemma~\ref{lem:formes_Pade}, we obtain that
$\Gamma(z)$ is non-singular and that
$
{\rm ord}_{z=1}\Delta(z)\leq \varepsilon_0 n+c_1$,
where 
$\Delta(z):=\det\Gamma(z)$, 
for some positive real number $c_1$ independent of $n$ and
$\varepsilon$. 
Let
$
c_2:={\rm ord}_{z=1}\Delta(z)$
and $
\nu:=i_3+c_2$; 
Then $c_2+\mu\leq\nu$, and, after enlarging $c_1$ if necessary,
$
\nu\leq\varepsilon_0 n+c_1$.

Set
\[
\mathbf L(z,\X):=
\prescript{\rm t}{}{
\bigl(
R_0(z,\X),\ldots,R_{\mu-1}(z,\X),
\ell_\mu(z,\X),\ldots,\ell_{i_3-1}(z,\X)
\bigr)}.
\]
By construction,
\[
\mathbf L(z,\X)
=
\Gamma(z)\,
\prescript{\rm t}{}{(X_0,\ldots,X_{i_3-1})}.
\]
Hence, by the comatrix formula,
\[
\Delta(z)
\prescript{\rm t}{}{(X_0,\ldots,X_{i_3-1})}
=
\prescript{\rm t}{}{\operatorname{Com}}(\Gamma(z))\,\mathbf L(z,\X).
\]

Recall that $\partial R_k=R_{k+1}$, that $\partial\ell_i$ is a
$\mathbb Q[z]$-linear combination of
$\ell_\mu,\ldots,\ell_{i+1}$ for $i<i_1$, and that
$\partial\ell_i=\ell_{i+1}$ for $i\geq i_1$. Applying
$\partial^{c_2}$ to the preceding identity, we obtain, for every
$j\in\{0,\ldots,i_3-1\}$,
\begin{multline*}
\Delta^{(c_2)}(z)X_j
+
\sum_{h=0}^{c_2-1}
\binom{c_2}{h}
\Delta^{(h)}(z)X_{j+c_2-h}\\
=
\sum_{i=0}^{\mu+c_2-1}q_{i,j}(z)R_i(z,\X)
+
\sum_{i=\mu}^{\nu-1}r_{i,j}(z)\ell_i(z,\X),
\end{multline*}
for suitable polynomials $q_{i,j},r_{i,j}\in\mathbb Q[z]$. 
Specializing at $z=1$ and using
that $\Delta^{(h)}(1)=0$, for $h<c_2$, and 
$
\Delta^{(c_2)}(1)\neq0$, 
we obtain
\[
X_j
=
\sum_{i=0}^{\mu+c_2-1}\lambda_{i,j}R_i(1,\X)
+
\sum_{i=\mu}^{\nu-1}\gamma_{i,j}\ell_i(1,\X),
\qquad
0\leq j<i_3,
\]
for suitable $\lambda_{i,j},\gamma_{i,j}\in\mathbb Q$. Thus
$X_0,\ldots,X_{i_3-1}$ belong to the $\mathbb Q$-vector space spanned
by
\begin{equation}\label{eq:formes_R_et_A}
R_0(1,\X),\ldots,R_{\mu+c_2-1}(1,\X),
\ell_\mu(1,\X),\ldots,\ell_{\nu-1}(1,\X).
\end{equation}

We now prove by induction on $j$ that $X_j$ belongs to this space for
every $j\leq\nu-\eta-1$. This has already been established for
$j<i_3$. Let
$
i_3\leq j\leq\nu-\eta-1$ 
and assume that the result holds for all smaller indices. Since
$j+\eta\geq i_2$, the coefficient $\lambda$ of $X_j$ in
$\ell_{j+\eta}(1,\X)$ is non-zero. Moreover,
$\ell_{j+\eta}(1,\X)$ is supported on $X_0,\ldots,X_j$, and therefore
\[
\ell_{j+\eta}(1,\X)
=
\lambda X_j+\sum_{k=0}^{j-1}\lambda_kX_k
\]
for suitable $\lambda_k\in\mathbb Q$. Since $j+\eta\leq\nu-1$, the
induction hypothesis shows that $X_j$ belongs to the space spanned
by~\eqref{eq:formes_R_et_A}. 
Consequently, the forms in~\eqref{eq:formes_R_et_A} span the space of
linear forms in $X_0,\ldots,X_{\nu-\eta-1}$. Since the forms
$
\ell_\mu(1,\X),\ldots,\ell_{\nu-1}(1,\X)$ 
are linearly independent, their span has dimension $\nu-\mu$. The
quotient by this span therefore has dimension $\mu-\eta$. We may
thus select $\mu$ forms
$$
L_1(\X),\ldots,L_\mu(\X)$$ 
among
$
R_0(1,\X),\ldots,R_{\mu+c_2-1}(1,\X)$ 
whose images span this quotient. These forms satisfy~(c), while
(a) and~(b) follow from their construction.
\end{proof}

Let
\[
M_n:=(a_{n;i,j})_{1\leq i,j\leq\nu}
\]
be the matrix defined by
\[
a_{n;i,j}:=
\begin{cases}
\text{the coefficient of $X_{j-1}$ in $L_i(\X)$},
& 1\leq i\leq\mu,\\
\text{the coefficient of $X_{j-1}$ in
$\ell_{i-1}(1,\X)$},
& \mu<i\leq\nu.
\end{cases}
\]
Then $(M_n)_{n\geq n_1}$ is our sequence of preparation matrices.
Unlike those considered in Section~\ref{sec:prepa_matrix}, the
matrices $M_n$ may not be invertible, as occurs when $\eta \geq 1$.

\subsubsection{The putative forms}

We now assume, for contradiction, that every differential equation
satisfied by $f$ has a singularity at $1$. In particular,
$a_{i_1}(1)=0$, and hence $\eta\geq1$.

Let $\kappa :=\mu-\eta$. Let $c_3$ be a positive real number, independent of $\varepsilon$, such
that
\begin{equation}\label{eq:majV_outstandingl_extra}
c_0(\varepsilon n+c_1+\eta-1)^\delta
\leq
(\varepsilon n)^{c_3}
\end{equation}
for every $n\geq2n_1\geq2\varepsilon^{-1}$. Recall that $i_0$ is the
least non-negative integer such that
$
f^{(i_0)}(1)\neq0$. 
Since $\varepsilon \kappa \leq \varepsilon\mu<1$ and $\eta\geq1$, we have
$
\eta-\varepsilon\kappa >0$. 
We now choose $\varepsilon>0$ sufficiently small that
$
\eta-\varepsilon\kappa-(c_3+1)\varepsilon>0$. 
Using $k!=\mathcal O(k^{k+1/2})$, we may then choose an integer
$n>2n_1$ such that
\begin{equation}\label{eq:cond_n_large}
(\varepsilon n+c_1-1)!\kappa\rho_1
\rho_2^{\kappa-1}
n^{-(\eta-\varepsilon\kappa)n} 
(\varepsilon n)^{c_3(\varepsilon n+c_1-\kappa)}
<
\bigl|f^{(i_0)}(1)\bigr|.
\end{equation}
We fix such an integer $n$ and consider the preparation matrix
associated with it.

The forms
$$
\ell_{\nu+i}(1,\X), \quad 0\leq i\leq\eta-1,$$ 
are supported on the variables $X_0,\ldots,X_{\nu-1}$. They constitute our set of putative forms.

\subsubsection{End of the proof of Theorem~\ref{thm:Andre_surQ}}

 Arguing as at the end of the proof of
Lemma~\ref{lem:sing_good_forms}, we obtain that
$X_{\nu-\eta},\ldots,X_{\nu-1}$ belong to the $\mathbb Q$-vector
space spanned by
$$
L_1(\X),\ldots,L_\mu(\X),
\ell_\mu(1,\X),\ldots,\ell_{\nu+\eta-1}(1,\X).
$$ 
Consequently, after reindexing the forms $L_i$ if necessary, the $\nu$ forms
\[
L_1(\X),\ldots,L_{\kappa}(\X),
\ell_\mu(1,\X),\ldots,\ell_{\nu+\eta-1}(1,\X)
\]
form a basis of the space of linear forms in
$X_0,\ldots,X_{\nu-1}$, where we recall that $\kappa=\mu-\eta$.

Let $M_n'$ denote the coefficient matrix of  these forms with respect to the basis
$X_0,\ldots,X_{\nu-1}$. It is invertible and it will be our forbidden matrix.  Set
$$
\bw:=
\prescript{\rm t}{}{
\bigl(f(1),f'(1),\ldots,f^{(\nu-1)}(1)\bigr)}.$$ 
Since
$
\ell_i(1,\f(1))=0$, for every $i\geq\mu$, 
we have
\begin{equation}\label{eq:singforb}
M_n'\bw
=
\prescript{\rm t}{}{
\bigl(e_{n,1},\ldots,e_{n,\kappa},0,\ldots,0\bigr)}=:\boldsymbol\varepsilon_n,
\end{equation}
and 
$
\|\boldsymbol\varepsilon_n\|
\leq
\rho_1n^{-(\mu-1-\varepsilon)n}$, by Lemma~\ref{lem:sing_good_forms}. 

Let $C_n$ denote the maximum of the absolute values of the entries in
the first $\kappa$ columns of $(M_n')^{-1}$. Arguing as in the proof
of Theorem~\ref{thm:SS}, we obtain
\begin{equation}\label{eq:majo_w}
\|  (M'_n)^{-1}(\boldsymbol\varepsilon_n)\|_\infty
\leq
\kappa C_n\rho_1
n^{-(\mu-1-\varepsilon)n}.
\end{equation} 
By \eqref{eq:majV_outstandingl_extra}, 
Lemma~\ref{lem:rel_func}, and the fact that $\nu\leq\varepsilon_0 n+c_1\leq\varepsilon n+c_1$, we have
\[
\house{\ell_i(1,\X)}
\leq
(\varepsilon n)^{c_3}
\qquad
(i\leq\nu+\eta-1).
\]
Moreover, Lemma~\ref{lem:sing_good_forms} gives
$
\house{L_i}
\leq
\rho_2n^{(1+\varepsilon)n}$. 
Since $\det M_n'$ is a non-zero integer, we have
$|\det M_n'|\geq1$. The comatrix formula therefore yields
\begin{equation}\label{eq:majo_Cn_2}
C_n
\leq
(\nu-1)!
\bigl(\rho_2n^{(1+\varepsilon)n}\bigr)^{\kappa-1}
(\varepsilon n)^{c_3(\nu-\kappa)}.
\end{equation}

Combining \eqref{eq:majo_w} and \eqref{eq:majo_Cn_2}, and using
$\nu\leq\varepsilon n+c_1$, we obtain
\begin{align*}
\| (M'_n)^{-1}\boldsymbol\varepsilon_n \|_\infty
&\leq
\kappa\rho_1n^{-(\mu-1-\varepsilon)n}
(\nu-1)!
\bigl(\rho_2n^{(1+\varepsilon)n}\bigr)^{\kappa-1}
(\varepsilon n)^{c_3(\nu-\kappa)}
\\
&\leq
(\varepsilon n+c_1-1)!\kappa\rho_1
\rho_2^{\kappa-1}
n^{-(\eta-\varepsilon\kappa)n}
(\varepsilon n)^{c_3(\varepsilon n+c_1-\kappa)}
\\
&<
\bigl|f^{(i_0)}(1)\bigr|,
\end{align*}
where the last inequality follows from
\eqref{eq:cond_n_large}. Since $i_0\leq\nu-1$, 
the definition of $\|\bw\|_\infty$ contradicts \eqref{eq:singforb}, showing that $M'_n$ is a forbidden matrix. 
This ends the proof. 
\qed

\subsection{Proof of the general case}

We could prove Theorem~\ref{thm:Andre} along the same lines as
Theorem~\ref{thm:Andre_surQ}. However, more than one additional linear
form would then be needed to obtain a contradiction, forcing us to
consider monomials in $f$ and its derivatives. Since this would lead to
rather cumbersome notation, we prefer to derive it from Theorem~\ref{thm:Andre_surQ}, using
elementary properties of differential operators.

We first recall the following classical lemma, for which we provide a proof for completeness.

\begin{lem}\label{lem:common_multiple}
Let $\alpha\in\Q$ and let
$L_1,\ldots,L_h\in\Q(z)\langle\partial\rangle$ be differential
operators having $\alpha$ as an ordinary point. Then there exists a
common left multiple of $L_1,\ldots,L_h$ in
$\Q(z)\langle\partial\rangle$ having $\alpha$ as an ordinary point.
\end{lem}

\begin{proof}
It is enough to consider two operators $L_1$ and $L_2$. Let $M$ denote
their monic least common left multiple. In a common differential
extension, the solution space of $M$ is the sum of the solution spaces
of $L_1$ and $L_2$. Indeed, if $G$ denotes their greatest common right
divisor, then
\[
\operatorname{ord}_\partial M
=
\operatorname{ord}_\partial L_1
+
\operatorname{ord}_\partial L_2
-
\operatorname{ord}_\partial G,
\]
while the intersection of the solution spaces of $L_1$ and $L_2$ is
the solution space of $G$.

Since $\alpha$ is an ordinary point of $L_1$ and $L_2$, all their local
solutions are holomorphic at $\alpha$. Hence all local solutions of $M$
are holomorphic at $\alpha$. Thus $\alpha$ is either an ordinary point
or an apparent singularity of $M$. By
\cite[Appendix, Prop.\ 7]{ABH06}, whose proof applies verbatim over the
algebraically closed field $\Q$, the operator $M$ admits a left
multiple $\widetilde M\in\Q(z)\langle\partial\rangle$ having $\alpha$
as an ordinary point. Since $\widetilde M$ is also a left multiple of
$L_1$ and $L_2$, the result follows. The general case follows by
induction.
\end{proof}

\begin{proof}[Proof of Theorem~\ref{thm:Andre}]
We first establish the following claim: for every $\alpha\in\C^\times$,
every $E$-function is annihilated by a differential operator in
$\Q(z)\langle\partial\rangle$ having $\alpha$ as an ordinary point.

Let $f$ be an $E$-function. If $\alpha\notin\Q$, then $\alpha$ is an
ordinary point of every differential operator in
$\Q(z)\langle\partial\rangle$, since the finite singularities of such
an operator are algebraic. The claim therefore holds in this case.

Suppose now that $\alpha\in\Q^\times$. Replacing $f(z)$ by
$f(\alpha z)$, it is enough to treat the case $\alpha=1$. Since the
coefficients of $f$ belong to a number field, we may write
\[
f=\beta_1f_1+\cdots+\beta_hf_h,
\]
where $\beta_1,\ldots,\beta_h\in\Q$ and
$f_1,\ldots,f_h\in\QQ[[z]]$ are $E$-functions (see \cite[Lem.\ 19]{AF24}). By
Theorem~\ref{thm:Andre_surQ}, for every $i\in\{1,\ldots,h\}$ there
exists an operator
$
L_i\in\QQ(z)\langle\partial\rangle$
such that $L_i(f_i)=0$ and $1$ is an ordinary point of $L_i$.
By Lemma~\ref{lem:common_multiple}, the operators
$L_1,\ldots,L_h$ admit a common left multiple $L$ in
$\Q(z)\langle\partial\rangle$ having $1$ as an ordinary point. Since
$L(f_i)=0$ for every $i$, we obtain
\[
L(f)
=
\beta_1L(f_1)+\cdots+\beta_hL(f_h)
=
0.
\]
After undoing the change of variable $z\mapsto\alpha z$, this proves
the claim.

Let now $M\in\Q(z)\langle\partial\rangle$ be the monic minimal
differential operator annihilating $f$. We claim that every non-zero
finite singularity of $M$ is apparent. Indeed, let
$\alpha\in\C^\times$. By the first part of the proof, there exists an
operator $L_\alpha\in\Q(z)\langle\partial\rangle$ annihilating $f$ and
having $\alpha$ as an ordinary point. By the minimality of $M$, there
exists $Q_\alpha\in\Q(z)\langle\partial\rangle$ such that
$
L_\alpha=Q_\alpha M$.
Consequently, every local solution of $M$ is also a local solution of
$L_\alpha$, and is therefore holomorphic at $\alpha$. Hence, if
$\alpha$ is a singularity of $M$, it is an apparent singularity.

We may now apply the desingularization theorem
\cite[Appendix, Thm.\ 6]{ABH06} to $M$. Its proof applies verbatim over
the algebraically closed field $\Q$ and yields a monic left multiple
$\widetilde M\in\Q(z)\langle\partial\rangle$ of $M$ such that every
ordinary point and every apparent singularity of $M$ is an ordinary
point of $\widetilde M$. Since all the non-zero finite singularities of
$M$ are apparent, the operator $\widetilde M$ has no non-zero finite
singularity. Finally, since $M(f)=0$ and $\widetilde M$ is a left
multiple of $M$, we have $\widetilde M(f)=0$. This proves the theorem.
\end{proof}

\begin{proof}[Proof of Andr\'e's theorem~{\rm (A1)}]
Let $f$ be an $E$-function. As observed in the proof of
Theorem~\ref{thm:Andre}, any point $\alpha\in\C^\times$ is either an
ordinary point or an apparent singularity of the minimal differential
operator annihilating $f$. Therefore, the only possible non-apparent
singularities of this operator are $0$ and $\infty$.
\end{proof}


\subsection{Proof of Corollary~\ref{thm:measure_sans_system}}
\label{sec:measures_E}

We are now ready to deduce
Corollary~\ref{thm:measure_sans_system} from
Theorems~\ref{thm:measure_system} and~\ref{thm:Andre}.

\begin{proof}[Proof of Corollary~\ref{thm:measure_sans_system}]
Every element of $\mathbf E$ can be written as the value at $z=1$ of an
$E$-function. For each $i\in\{1,\ldots,r\}$, let $f_i$ be an
$E$-function such that $\xi_i=f_i(1)$. By
Theorem~\ref{thm:Andre}, the function $f_i$ satisfies a differential
equation
\begin{equation}\label{eq:goodE}
a_{i,0}(z)f_i(z)+a_{i,1}(z)f_i'(z)+\cdots
+a_{i,m_i}(z)f_i^{(m_i)}(z)=0
\end{equation}
having no non-zero finite singularity. In particular,
$a_{i,m_i}(1)\neq 0$. 
The column vector
$\prescript{\rm t}{}{(f_i(z),f_i'(z),\ldots,f_i^{(m_i-1)}(z))}$
is a solution of the linear differential system
$Y'(z)=A_i(z)Y(z)$, where
\begin{equation}\label{eq:Ai}
A_i(z):=
\begin{pmatrix}
0 & 1 & 0 & \cdots & 0 \\
\vdots & \ddots & \ddots & & \vdots \\
\vdots & & \ddots & \ddots & 0 \\
0 & \cdots & \cdots & 0 & 1 \\
-\dfrac{a_{i,0}(z)}{a_{i,m_i}(z)}
& \cdots
& \cdots
& \cdots
& -\dfrac{a_{i,m_i-1}(z)}{a_{i,m_i}(z)}
\end{pmatrix}.
\end{equation}
The functions $f_i^{(k)}$, with $1\leq i\leq r$ and
$0\leq k\leq m_i-1$, are therefore the coordinates of a vector solution
of the differential system $Y'(z)=A(z)Y(z)$, where $A(z)$ is the
block-diagonal matrix
\[
A(z)
:=
A_1(z)\oplus\cdots\oplus A_r(z)
=
\begin{pmatrix}
A_1(z) & & \\
& \ddots & \\
& & A_r(z)
\end{pmatrix}.
\]
The point $1$ is regular with respect to this system. The result now
follows directly from Theorem~\ref{thm:measure_system}. 
Moreover, $t$ is equal to the transcendence degree over $\C(z)$ of the
functions $f_i^{(k)}$, with $1\leq i\leq r$ and
$0\leq k\leq m_i-1$, or equivalently, to the transcendence degree over
$\C(z)$ of the differential field generated by $f_1,\ldots,f_r$.
\end{proof}

\subsubsection{Remarks on Corollary~\ref{thm:measure_sans_system}}

Let us end this section with a few comments on
Corollary~\ref{thm:measure_sans_system}.

Recall that a real transcendental number $\xi$ is called a
\emph{$U_d$-number} if, for every $w\geq 1$, there exist non-zero
polynomials $P\in\Z[X]$ of degree at most $d$ and height at most 
$H\geq 2$ such that
\[
0<|P(\xi)|\leq H^{-w}.
\] 
A real transcendental number is
called a \emph{$U$-number} if it is a $U_d$-number for some positive
integer $d$. In particular, the $U_1$-numbers are exactly the Liouville
numbers.

The question of whether some values of $E$-functions may be
$U$-numbers, and in particular Liouville numbers, is a longstanding
problem with roots in the pioneering work of Siegel~\cite{Si29}. It was fully
resolved only recently by Fischler and Rivoal~\cite{FR23}, who used a
desingularization procedure due to Beukers together with a
Galois-theoretic argument. Their proof relies on Beukers' lifting
theorem and on Andr\'e's result~(A1). As a direct consequence of
Corollary~\ref{thm:measure_sans_system}, we obtain a new proof of their result, one that does not invoke Beukers' lifting theorem.

\begin{coro}\label{coro:Unumbers}
No value of an $E$-function at an algebraic point is a $U$-number.
In particular, no such value is a Liouville number.
\end{coro}

\begin{proof}
Let $\xi\in\mathbf E$. If $\xi$ is not real or if it is algebraic,
there is nothing to prove. Assume therefore that $\xi$ is real and
transcendental. Let $d\geq 1$ be fixed. By
Corollary~\ref{thm:measure_sans_system}, there exist positive constants
$A$ and $B$, depending only on $\xi$ and $d$, such that, for every
non-zero polynomial $P\in\Z[X]$ of degree at most $d$ and height at most $H$,
\[
|P(\xi)|\geq A H^{-B}.
\]
It follows immediately that $\xi$ is not a $U_d$-number. Since $d$ is
arbitrary, $\xi$ is not a $U$-number.
\end{proof}

A general form of Liouville's inequality may be stated as follows. Given
algebraic numbers $\alpha_1,\ldots,\alpha_m$ and a polynomial
$P\in\Z[X_1,\ldots,X_m]$ of (total) degree at most $d$ height at most $H$, one has
\[
\text{either } P(\alpha_1,\ldots,\alpha_m)=0,
\qquad
\text{or }
|P(\alpha_1,\ldots,\alpha_m)|
\geq C_1H^{-C_2},
\]
where $C_1$ depends on $d$ and
$\alpha_1,\ldots,\alpha_m$, whereas $C_2$ depends only on
$\alpha_1,\ldots,\alpha_m$; see, for instance,
\cite[p.~84--85]{Wa_Liv}. In this respect,
Corollary~\ref{thm:measure_sans_system} may be viewed as a
Liouville-type inequality for elements of $\mathbf E$. 
This estimate allows one to work with elements of $\mathbf E$ in much
the same way as with algebraic numbers, though the bounds are far less
sharp. In the case of values of the Bessel function, this was already
observed by Siegel, just after stating his main result~\cite{Si29}
(see also \cite[pp.~29--30]{Si14}).

Just as Liouville's inequality implies the transcendence of
$\sum_{n\geq 0}10^{-n!}$, Corollary~\ref{thm:measure_sans_system},
together with the effective estimates discussed in
Section~\ref{sec:effectivity}, implies that a number of the form
\[
\sum_{n=0}^{\infty}e^{-u_n}
\]
is transcendental over $\mathbf E$, provided that the strictly
increasing sequence of positive integers $(u_n)_{n\geq 0}$ grows
sufficiently rapidly. More precisely, it is enough to assume that, for
every real number $A>0$,
\[
\lim_{n\to\infty}
\frac{\log u_{n+1}}{u_n^A\log u_n}
=
+\infty.
\]
For example, one may take $u_n:=\prescript{n}{}{n}$, where
\[
\prescript{n}{}{n}
\;:=\;
\underbrace{n^{n^{{\cdot}^{{\cdot}^{{\cdot}^{n}}}}}}_{n},
\]
denotes the power tower of $n$ copies of $n$.

\appendix

\renewcommand{\thesection}{\Roman{section}}

\section{Specialization of differential ideals}
\label{AppendiceA}

The main purpose of this appendix is to provide an elementary proof of Andr\'e's result~{\rm (A2)}
stated in the Introduction. It will be deduced from the more general
Proposition~\ref{prop:specialization_ideal} below, which is itself a
consequence of \cite[Thm.~1.6.1(3)]{An14}.

Recall that, for a field $k$, an ideal
$\I\subset k[X_1,\ldots,X_m]$ is said to be
\emph{absolutely prime} if, for every field extension
$K$ of $k$, the ideal generated by $\I$ in
$K[X_1,\ldots,X_m]$ remains prime. This is equivalent to the fact that
the fraction field of $k[X_1,\ldots,X_m]/\I$ is a regular extension of
$k$ (see, for instance, \cite[Thm.~39]{ZS}). Moreover, when
$k$ has characteristic zero, this is equivalent to the fact that
$k$ is algebraically closed in the fraction field of
$k[X_1,\ldots,X_m]/\I$.

Let $\Q[z]_{(z-\alpha)}$ denote the localization of $\Q[z]$
at the prime ideal $(z-\alpha)$, that is, the ring of rational
functions that are regular at~$\alpha$. An element
$P\in\Q[z]_{(z-\alpha)}[X_1,\ldots,X_m]$ can naturally be specialized
at $\alpha$ by setting
$
{\rm ev}_\alpha(P):=P(\alpha,X_1,\ldots,X_m)
\in\Q[X_1,\ldots,X_m]$. 
For an ideal $\I\subset\Q(z)[X_1,\ldots,X_m]$, we set 
$$
{\rm ev}_\alpha(\I):=
{\rm ev}_\alpha\bigl(\I\cap\Q[z]_{(z-\alpha)}[X_1,\ldots,X_m]\bigr).
$$

\begin{rem}
We have
$
{\rm ev}_\alpha(\I)
={\rm ev}_\alpha\bigl(\I\cap\Q[z][X_1,\ldots,X_m]\bigr).
$
Indeed, if
$
P=\frac{Q}{d(z)}\in
\Q[z]_{(z-\alpha)}[X_1,\ldots,X_m],
$
where $Q\in\Q[z][X_1,\ldots,X_m]$ and $d(\alpha)\neq0$, then
$
{\rm ev}_\alpha(P)
=
{\rm ev}_\alpha\!\left(\frac{Q}{d(\alpha)}\right).
$
\end{rem}

\begin{prop}\label{prop:specialization_ideal}
Let $D$ be a derivation on
$\Q(z)[\X]=\Q(z)[X_1,\ldots,X_m]$ acting on $\Q(z)$ as the usual
derivation and on
$\X:=\prescript{\rm t}{}{(X_1,\ldots,X_m)}$ by
\begin{equation}\label{eq:D}
D\X=A(z)\X
\end{equation}
for some matrix $A(z)\in\mathcal M_m(\Q(z))$. Let
$\I\subset\Q(z)[\X]$ be a $D$-invariant absolutely prime ideal, and let
$\alpha\in\Q$ be a point that is not a pole of $A(z)$. Then
${\rm ev}_\alpha(\I)$ is prime.
\end{prop}

By a classical result of Krull~\cite{Kr}, the specialization of an absolutely prime ideal $\I$ of 
$\Q(z)[X_1,\ldots,X_m]$ remains prime for all but finitely many $\alpha$. Proposition~\ref{prop:specialization_ideal} strengthens this by showing
that, when $\I$ is $D$-invariant, this finite exceptional set is actually
contained in the set of singular points of the differential
system~\eqref{eq:D}.
 
Before proving Proposition~\ref{prop:specialization_ideal}, we show how it
implies Andr\'e's result~{\rm (A2)}. Let $\mathfrak I$ be the ideal
appearing in~{\rm (A2)}, and let 
$
\mathfrak I_{\Q(z)}
:=
\mathfrak I\,\Q(z)[X_1,\ldots,X_m]$ be its extension to $\Q(z)[X_1,\ldots,X_m]$. 
 This ideal is invariant
under the derivation $D$ defined by $D\X=A(z)\X$. 
Moreover, the field
$
\Q(z)(f_1,\ldots,f_m)$ 
is a regular extension of $\Q(z)$. Indeed, since $E$-functions are entire,
$\Q(z)(f_1,\ldots,f_m)$ is a field of single-valued meromorphic
functions on $\mathbb C$. Hence $\Q(z)$ is algebraically closed in 
$\Q(z)(f_1,\ldots,f_m)$ and the 
extension is regular, so that $\mathfrak I_{\Q(z)}$ is absolutely prime. 
Proposition~\ref{prop:specialization_ideal} now shows that
${\rm ev}_\alpha(\mathfrak I_{\Q(z)})$ is prime. Since
${\rm ev}_\alpha(\mathfrak I_{\Q(z)})
=
{\rm ev}_\alpha(\mathfrak I)$,
the latter ideal is prime.

\subsection{Proof of Proposition~\ref{prop:specialization_ideal}}

We now turn to the proof of Proposition~\ref{prop:specialization_ideal}. By
replacing $z$ with $z+\alpha$, it suffices to prove the result for
$\alpha=0$. The proof proceeds in two steps. We first establish the
result when the derivation is the trivial extension of $\frac{d}{dz}$. We then reduce
the general case to this special case by using a gauge transformation
associated with a fundamental matrix of solutions. Since this matrix
has entries in $\Q[[z]]$, the corresponding change of variables is
naturally carried out over the Laurent series field $\Q((z))$.

We identify $\Q[[z]][\X]$ with a subring of $\Q[\X][[z]]$ in a natural way. Thus every element
$P\in\Q[[z]][\X]$ admits a unique expansion
\[
P=\sum_{n\ge0}P_n(\X)z^n,
\qquad P_n(\X)\in\Q[\X],
\]
and we define the specialization morphism
\[
{\rm ev}_0:\Q[[z]][\X]\longrightarrow\Q[\X]
\]
by
$
{\rm ev}_0(P):=P_0(\X)$. 
For an ideal $\I\subset\Q((z))[\X]$, we set
\[
{\rm ev}_0(\I):=
{\rm ev}_0\bigl(\I\cap\Q[[z]][\X]\bigr).
\]
We extend the derivation $\frac{d}{dz}$ from $\Q((z))$ to $\Q((z))[\X]$ by defining
$$
\frac{d}{dz} X_i=0,\quad 1\leq i \leq m.
$$

\begin{lem}\label{lem:specialization_derivative}
Let $\I\subset\Q((z))[\X]$ be a
$\frac{d}{dz}$-invariant ideal. Then
\begin{equation}\label{eq:specialization_derivative}
{\rm ev}_0(\I)=\I\cap\Q[\X].
\end{equation}
In particular, if $\I$ is prime, then so is ${\rm ev}_0(\I)$.
\end{lem}

\begin{proof}
The final assertion follows immediately from
\eqref{eq:specialization_derivative}, since the contraction of a prime
ideal is prime. It remains to prove the equality.

The inclusion
$
\I\cap\Q[\X]\subseteq{\rm ev}_0(\I)$ 
is immediate.  
Conversely, let
$
Q(\X)\in{\rm ev}_0(\I)$. 
By definition, there exists
$
R(z,\X)\in\I\cap\Q[[z]][\X]$ 
such that
$
Q(\X)={\rm ev}_0(R)$.  
Let $d$ be the total degree of $R$ in $\X$, and consider a
$\Q((z))$-basis
\[
P_1(z,\X),\ldots,P_r(z,\X)
\]
of the vector space of elements of $\I$ having total degree at most
$d$ in $\X$. After multiplying the $P_i$ by suitable powers of $z$, we
may assume that
$
P_1,\ldots,P_r\in\Q[[z]][\X]$. 

The proof of \cite[Lem.~3.1]{Be06} carries over \textit{mutatis mutandis}
to linear forms with coefficients in $\Q[[z]]$.  As a consequence, we can choose such a basis so that
$
{\rm ev}_0(P_1),\ldots,{\rm ev}_0(P_r)$ 
are linearly independent over $\Q$.

Since $R$ has total degree at most $d$ in $\X$, there exist
$g_1(z),\ldots,g_r(z)\in\Q((z))$ such that
\[
R(z,\X)
=
g_1(z)P_1(z,\X)+\cdots+g_r(z)P_r(z,\X).
\]
The linear independence of
${\rm ev}_0(P_1),\ldots,{\rm ev}_0(P_r)$ implies that all the $g_i$
belong to $\Q[[z]]$. Indeed, otherwise, multiplying this equality by
the smallest power of $z$ that removes all the poles and then applying
${\rm ev}_0$ would yield a non-trivial linear relation among
${\rm ev}_0(P_1),\ldots,{\rm ev}_0(P_r)$. Thus
\begin{equation}\label{eq:Q}
Q(\X)
=
\sum_{i=1}^r g_i(0)\,{\rm ev}_0(P_i).
\end{equation}

Let $\bP(z,\X)$ denote the column vector with entries
$P_1(z,\X),\ldots,P_r(z,\X)$. Since $\I$ is invariant under
$\frac{d}{dz}$ and this derivation does not increase the
degree in $\X$, there exist matrices
$
B_0(z)={\rm I}_r,B_1(z),B_2(z),\ldots$ 
with entries in $\Q((z))$ such that
\[
\bP^{(k)}(z,\X)=B_k(z)\bP(z,\X),
\qquad k\geq0.
\]
The linear independence of the entries of ${\rm ev}_0(\bP)$ implies
that every $B_k(z)$ has entries in $\Q[[z]]$.  
Expanding the entries of $\bP$ into their formal Taylor series, we obtain
\[
\bP(z,\X)
=
\sum_{k=0}^{\infty}\frac{z^k}{k!}\bP^{(k)}(0,\X)
=
\left(
\sum_{k=0}^{\infty}\frac{z^k}{k!}B_k(0)
\right)
{\rm ev}_0(\bP).
\]
The matrix
\[
C(z)
:=
\sum_{k=0}^{\infty}\frac{z^k}{k!}B_k(0)
\]
has entries in $\Q[[z]]$ and satisfies $C(0)={\rm I}_r$. Hence
$C(z)\in\GL_r(\Q[[z]])$ and
\[
{\rm ev}_0(\bP)=C(z)^{-1}\bP(z,\X).
\]
It follows that
$
{\rm ev}_0(P_1),\ldots,{\rm ev}_0(P_r)\in\I$. 
Equation~\eqref{eq:Q} therefore implies that $Q\in\I$. Thus
$
{\rm ev}_0(\I)\subseteq\I\cap\Q[\X]$, 
which completes the proof.
\end{proof}

We now show that the general case reduces to
Lemma~\ref{lem:specialization_derivative}. The derivation $D$
extends to $\Q((z))[\X]$ by setting
$
Dg:=g'(z)$  for every $g\in\Q((z))$, and 
$D\X:=A(z)\X$.   
Consider the differential system
$
Y'(z)=A(z)Y(z)$. 
Since $A(z)$ has no pole at $0$, there exists a 
unique  \emph{fundamental matrix of solutions} 
$
F(z)\in\GL_m(\Q[[z]])$ 
such that
\[
F'(z)=A(z)F(z)
\qquad\text{and}\qquad
F(0)={\rm I}_m.
\]
Define
\[
\begin{aligned}
\phi_F:\Q((z))[\X]&\longrightarrow\Q((z))[\X],\\
P(z,\X)&\longmapsto P(z,F(z)\X).
\end{aligned}
\]
This is an automorphism, with inverse $\phi_{F^{-1}}$, and it restricts
to an automorphism of $\Q[[z]][\X]$. 

\begin{lem}\label{lem:D-stable}
Let $\I\subset\Q((z))[\X]$ be a $D$-invariant ideal. Then
$\phi_F(\I)$ is $\frac{d}{dz}$-invariant.
\end{lem}

\begin{proof}
We claim that, for every $P\in\Q((z))[\X]$,
$$
\phi_F(DP)
=
\frac{d}{dz}\phi_F(P).$$ 
The two sides agree on $\Q((z))$ and satisfy the same Leibniz rule
relative to $\phi_F$. It is therefore enough to check the identity on
$X_1,\ldots,X_m$. In vector notation, we have
\[
\phi_F(D\X)
=
\phi_F(A(z)\X)
=
A(z)F(z)\X
=
F'(z)\X
=
\frac{d}{dz}(F(z)\X)
=
\frac{d}{dz}\phi_F(\X).
\]
Now let $P\in\phi_F(\I)$. Then $P=\phi_F(Q)$ for some $Q\in\I$, and
therefore
\[
\frac{d P}{d z}
=
\frac{d}{dz}\phi_F(Q)
=
\phi_F(DQ)
\in
\phi_F(\I),
\]
since $DQ\in\I$.
\end{proof}

We may now prove Proposition~\ref{prop:specialization_ideal}.

\begin{proof}[Proof of Proposition~\ref{prop:specialization_ideal}]
By assumption, $\I$ is absolutely prime. Hence the extended ideal
$
\widehat{\I}:=\I\,\Q((z))[\X]$ 
is prime. Since $\phi_F$ is an automorphism, the ideal
$\phi_F(\widehat{\I})\subset\Q((z))[\X]$ is also prime.
By Lemma~\ref{lem:D-stable}, it is
$\frac{d}{dz}$-invariant. By
Lemma~\ref{lem:specialization_derivative}, the specialization
${\rm ev}_0(\phi_F(\widehat{\I}))$ is therefore prime.

Now $\phi_F$ and $\phi_{F^{-1}}$ both preserve $\Q[[z]][\X]$.
Moreover, since $F(0)={\rm I}_m$, we have
$
{\rm ev}_0(\phi_F(P))={\rm ev}_0(P)$ 
for every $P\in\Q[[z]][\X]$. Consequently,
${\rm ev}_0(\phi_F(\widehat{\I}))
=
{\rm ev}_0(\widehat{\I})$, and hence ${\rm ev}_0(\widehat{\I})$ is prime. 

Finally, since $\widehat{\I}$ is the extension of $\I$ to $\Q((z))[\X]$, it follows that ${\rm ev}_0(\I)={\rm ev}_0(\widehat{\I})$.
Hence ${\rm ev}_0(\I)$ is prime, and the proof is complete.
\end{proof}

\section{Further consequences of Theorem~\ref{thm:Andre}}
\label{AppendiceB}

Although {\rm (A1)} is the only part of Andr\'e's result that has
been used so far in the transcendence theory of $E$-functions,
Andr\'e actually proved in~\cite{An1} a more precise description of the
singularities of the minimal differential equation of a non-polynomial
$E$-function: it is regular singular at $0$ and irregular at infinity,
with slopes in $\{0,1\}$, and it admits a full basis of solutions in
the class of Nilsson--Gevrey series of arithmetic type and order $-1$.
Obtaining this more precise result requires the full strength of the
theory of $G$-operators, including results of the Chudnovskys, Katz,
and Andr\'e himself. Here and in Theorem~\ref{thm:Andre_consequence}
below, the slopes refer to those of the local Newton polygon at infinity
of a differential equation as considered in~\cite{vdPS03}. 

Our approach does not seem to recover this full result. Nevertheless,
we show that Theorem~\ref{thm:Andre} easily implies the following
refinement of~{\rm (A1)}. In particular, the Siegel--Shidlovskii
method, as revisited in the present work, also shows that the minimal
differential equation satisfied by any $G$-function is regular singular
at the origin.

\begin{thm}\label{thm:Andre_consequence}
The minimal differential equation satisfied by a non-polynomial
$E$-function has no non-trivial singularities outside $\{0,\infty\}$.
Moreover,
\begin{itemize}
\item[-] it is regular singular at $0$;
\item[-] it is irregular at infinity, with slopes belonging to the
interval $[0,1]$.
\end{itemize}
Furthermore, the minimal differential equation satisfied by any
$G$-function is regular singular at $0$.
\end{thm}

Let us first recall some facts about differential operators. Set $\partial=d/dz$. We denote by $\mathcal F$
the Fourier--Laplace transform, that is, the automorphism of the Weyl
algebra $\mathbb C[z]\langle\partial\rangle$ defined by
$
\mathcal F(z):=-\partial
$
and
$
\mathcal F(\partial):=z.
$
It has order four, and its inverse, denoted by
$\overline{\mathcal F}:=\mathcal F^{-1}$, satisfies 
$
\overline{\mathcal F}(z)=\partial
$
and
$
\overline{\mathcal F}(\partial)=-z.
$

In the proof of Theorem~\ref{thm:Andre_consequence}, it will be more convenient to use
the Newton--Ramis polygon instead of the more classical local Newton polygons.
The Newton--Ramis polygon $\mathcal N(L)$ of a differential operator
\[
L:=a_m(z)\partial^m+a_{m-1}(z)\partial^{m-1}+\cdots+a_0(z)
\]
is defined as the convex hull of the horizontal half-lines whose right-hand
endpoints are
\[
(i,\val(a_i)-i)
\qquad\text{and}\qquad
(i,\deg(a_i)-i),
\]
for those $i$ such that $a_i\neq0$ (see an example below).
This convention differs from the one used for local Newton polygons of
differential equations in~\cite{vdPS03}:  the positive slopes
of $\mathcal N(L)$ are the finite positive slopes of the local Newton polygon
at $0$, whereas the opposites of the negative slopes of $\mathcal N(L)$ are
the finite positive slopes of the local Newton polygon at infinity, computed
after the change of variable $z\mapsto z^{-1}$; see
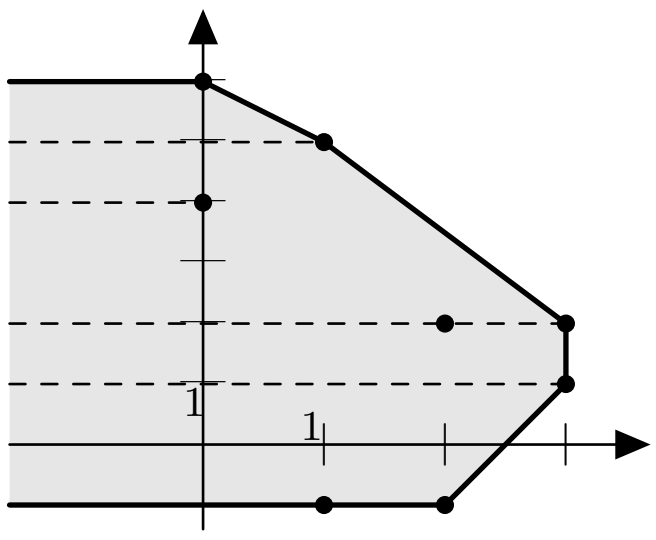
\cite[Rem.~p.~179]{Ramis}.

    \begin{figure}[h!]
	\centering
	\definecolor{rfsqsq}{rgb}{0,0,0}
	\begin{tikzpicture}[line cap=round,line join=round,>=triangle 45,x=0.8cm,y=0.4cm]
	\clip(-2,-1.5) rectangle (5,8);
	\fill[line width=1pt,color=rfsqsq,fill=rfsqsq,fill opacity=0.1] (-1.6,-1) -- (2,-1) -- (3,1) -- (3,2) -- (1,5) -- (0,6) -- (-1.6,6) -- cycle;
	\draw [line width=1pt] (-1.6,-1)-- (2,-1);
	\draw [line width=1pt] (2,-1)-- (3,1);
	\draw [line width=1pt] (3,1) -- (3,2);
	\draw [line width=1pt] (3,2) -- (1,5);
	\draw [line width=1pt] (1,5) -- (0,6);
	\draw [line width=1pt] (0,6) -- (-1.6,6);	
	\draw [->,line width=0.5pt] (-1.6,0) -- (3.7,0);
	\draw [->,line width=0.5pt] (0,-1.4) -- (0,7.2);
	\begin{scriptsize}
	\draw [fill=black] (1,-1) circle (1.5pt);
	\draw [line width=0.5pt, dashed] (-1.6,1)-- (3,1);
	\draw [line width=0.5pt, dashed] (-1.6,2)-- (3,2);
	\draw [line width=0.5pt, dashed] (-1.6,4)-- (0,4);
	\draw [line width=0.5pt, dashed] (-1.6,5)-- (1,5);
	\draw (-0.07,0.7) node {$1$};
	\draw (0.9,0.3) node {$1$};
	\draw [fill=black] (2,2) circle (1.5pt);
	\draw [fill=black] (2,-1) circle (1.5pt);
	\draw [fill=black] (3,1) circle (1.5pt);
	\draw [fill=black] (3,2) circle (1.5pt);
	\draw [fill=black] (1,5) circle (1.5pt);
	\draw [fill=black] (0,6) circle (1.5pt);
	\draw [fill=black] (0,4) circle (1.5pt);
	\draw (3,0) node {$\vert$};
	\draw (2,0) node {$\vert$};
	\draw (1,0) node {$\vert$};
	\draw (0,1) node {\textemdash};
	\draw (0,2) node {\textemdash};
	\draw (0,3) node {\textemdash};
	\draw (0,4) node {\textemdash};
	\draw (0,5) node {\textemdash};
	\draw (0,6) node {\textemdash};
	\end{scriptsize}
	\end{tikzpicture}
	\caption{\small\textit{The Newton-Ramis polygon of the differential operator $(3z^5-z^4)\partial^3+(z^4+z)\partial^2+(z^6+z^2+3z+1)\partial+(z^6+z^4)$.}}
\end{figure}

We shall use the following properties of the Newton--Ramis polygon.

\begin{itemize}
\item[(a)] If $L=L_1\cdots L_k$, then the set of finite slopes of
$\mathcal N(L)$ is the union of the sets of finite slopes of the
polygons $\mathcal N(L_i)$; see \cite[Lem.~3.45]{vdPS03}.

\item[(b)] The operator $L$ is regular singular at $0$, respectively
at infinity, if and only if $\mathcal N(L)$ has no positive,
respectively no negative, slope; see \cite[Ex.~3.47]{vdPS03}.

\item[(c)] Under the Fourier--Laplace transform, a slope $t$ is mapped
to
$
-\frac{t}{t+1};
$
see \cite[Sec.~5]{An1}. This transformation is understood as an
involution of $\mathbb R\cup\{\infty\}$, exchanging $-1$ and
$\infty$.

\item[(d)] A change of variable $z\mapsto z^\nu$, with
$\nu\in\mathbb Z\setminus\{0\}$, multiplies the slopes of the
Newton--Ramis polygon by $\nu$; see \cite[p.~180]{Ramis}. In
particular, the change of variable $z\mapsto z^{-1}$ changes the signs
of the slopes.

\item[(e)] Let $L$ be the least common left multiple of
$L_1,\ldots,L_k$. Then the set of finite slopes of $\mathcal N(L)$ is
the union of the sets of finite slopes of the polygons
$\mathcal N(L_i)$.

This statement is not explicitly formulated in \cite{vdPS03}, but it
follows from the formal decomposition theorem. Indeed, in a common
formal differential extension and after a finite ramification, the
solution space of an operator $M$ decomposes according to its
exponential parts: its solutions are finite sums of terms of the form
$f e(q)$, where
$
f\in
\mathbb C\bigl(\!\bigl(z^{1/\star}\bigr)\!\bigr)
\bigl[(z^a)_{a\in\mathbb C},\log z\bigr]
$
and $e(q)$ satisfies
$
z\frac{d}{dz}e(q)=q e(q),
$
with
$
q\in\bigcup_{r\geq1}z^{-1/r}\mathbb C[z^{-1/r}];
$
see \cite[p.~76]{vdPS03}. The non-negative slopes of $\mathcal N(M)$
are the opposites of the valuations of the exponential parts $q$; see
\cite[Thm.~3.54 and Rem.~3.55]{vdPS03}. In such a common formal
extension, the solution space of a least common left multiple is the
sum of the solution spaces of the operators involved. Its exponential
parts, and hence its non-negative slopes, are therefore exactly the
union of theirs. The corresponding assertion for the negative slopes
follows by applying the change of variable $z\mapsto z^{-1}$.

\item[(f)] Let
$
f(z):=\sum_{n\geq0}\frac{a_n}{n!}z^n$,
$
g(z):=\sum_{n\geq0}a_nz^n$,
and
$
h(z):=z^{-1}g(z^{-1}).
$
If $L$ annihilates $f$ and $\nu$ is sufficiently large, then
$\partial^\nu\mathcal F(L)$ annihilates $h$. Conversely, if $L$
annihilates $h$, then $\overline{\mathcal F}(L)$ annihilates $f$; see
\cite[Sec.~2]{An1}.
\end{itemize}

We shall also use the following elementary observation, for which we provide a proof for completeness.

\begin{lem}\label{lem:ramification_minimal}
Let $u$ be a non-zero holonomic power series and let $q$ be a positive integer.
If $L$ is the minimal operator annihilating $u(z)$, then the pullback
of $L$ under the change of variable $z\mapsto z^q$ is the minimal
operator annihilating $u(z^q)$. In particular, the slopes of the
minimal operator of $u(z^q)$ are obtained from those of $L$ by
multiplication by $q$.
\end{lem}

\begin{proof}
Set $\theta:=z\partial$. The order of the minimal differential
operator annihilating $u$ is the dimension over $\mathbb C(z)$ of the
vector space spanned by
$
u,\theta u,\theta^2u,\ldots.
$
Moreover,
$
\theta^j\bigl(u(z^q)\bigr)
=q^j(\theta^ju)(z^q)
$
for every $j\geq0$. 
We use the standard fact that $\mathbb C(z^{1/q})$ and
$\mathbb C((z))$ are linearly disjoint over $\mathbb C(z)$. Hence
elements of $\mathbb C((z))$ which are linearly independent over
$\mathbb C(z)$ remain linearly independent over
$\mathbb C(z^{1/q})$. 
Let $x:=z^q$. A linear relation over $\mathbb C(z)$ among
$
u(z^q),\theta u(z^q),\ldots,\theta^r u(z^q)
$
becomes a linear relation over $\mathbb C(x^{1/q})$ among
$
u(x),\theta u(x),\ldots,\theta^r u(x).
$
It follows that the minimal operators of $u(z)$ and $u(z^q)$ have the
same order. 
The pullback of $L$ under $z\mapsto z^q$ annihilates $u(z^q)$ and
has the same order as $L$. It is therefore minimal. The assertion
concerning the slopes follows from~(d).
\end{proof}

Our proof is divided into five steps.

\subsection*{Step 1}

\textit{Every $E$-function is annihilated by an operator whose
Newton--Ramis polygon has no infinite slope.}

\medskip

By definition, the Newton--Ramis polygon of an operator has an infinite slope if and only if $\val(a_m)<\deg(a_m)$. 
Theorem~\ref{thm:Andre} implies that every $E$-function is annihilated
by a differential operator with polynomial coefficients whose only
possible singularities are $0$ and infinity. Its leading coefficient
is therefore of the form $cz^\nu$, for some $c\neq0$ and some
non-negative integer $\nu$.
Hence
\[
\val(a_m)=\deg(a_m)=\nu,
\]
so the associated Newton--Ramis polygon has no infinite slope. 

\subsection*{Step 2}

\textit{The Newton--Ramis polygon of the minimal operator annihilating
a $G$-function has no slope equal to $1$.}

\medskip

Let
$
g(z):=\sum_{n\geq0}a_nz^n
$
and
$
f(z):=\sum_{n\geq0}\frac{a_n}{n!}z^n
$
be the associated $G$- and $E$-functions. Let $L$ be an operator
annihilating $f$ as in Step~1. Since the map
$
t\longmapsto-\frac{t}{t+1}
$
sends $\infty$ to $-1$, the Newton--Ramis polygon of
$\mathcal F(L)$ has no slope equal to $-1$. 
Choose $\nu$ such that
$
\Psi:=\partial^\nu\mathcal F(L)
$
annihilates $z^{-1}g(z^{-1})$, as in~(f). Since $\partial$ has no
finite non-zero slope, property~(a) shows that $\mathcal N(\Psi)$ still
has no slope equal to $-1$. 
Let $\Phi$ be the operator obtained from $\Psi$ by the change of
variable $z\mapsto z^{-1}$. Then $\Phi$ annihilates $zg(z)$ and, by
(d), $\mathcal N(\Phi)$ has no slope equal to $1$. Consequently, the
operator $\Phi z$, where $z$ denotes right multiplication by $z$,
annihilates $g$. Since the operator of multiplication by $z$ has no
finite non-zero slope, property~(a) shows that $\mathcal N(\Phi z)$ has
no slope equal to $1$. 
The minimal operator of $g$ is a right factor of every operator
annihilating $g$. Property~(a) therefore implies that its
Newton--Ramis polygon has no slope equal to $1$.

\subsection*{Step 3}

\textit{The Newton--Ramis polygon of the minimal operator annihilating
a $G$-function has no positive integral slope.}

\medskip

We argue by contradiction. Let $g$ be a $G$-function, let $L$ be its
minimal operator, and suppose that $\mathcal N(L)$ has an integral
slope $t\geq1$. Write
$
g(z):=\sum_{i=0}^{t-1}z^i g_i(z^t).
$
For every $i$, the series $z^ig_i(z^t)$ is the Hadamard product of
$g(z)$ and $z^i/(1-z^t)$, and is therefore a $G$-function. Moreover,
each $g_i$ is a $G$-function, since the class of $G$-functions is
stable under sections of their coefficient sequences. 
Let $L_i$ be the minimal operator annihilating
$z^ig_i(z^t)$, and let $M$ be the least common left multiple of
$L_0,\ldots,L_{t-1}$. Then $M$ annihilates $g$, so the minimal
operator $L$ is a right factor of $M$. By~(a), the slope $t$ of
$\mathcal N(L)$ is therefore a slope of $\mathcal N(M)$. By~(e), it
must be a slope of $\mathcal N(L_i)$ for at least one $i$. 
Multiplication of a solution by the monomial $z^i$ amounts to a
regular singular gauge transformation and does not change the slopes.
Hence the minimal operator of $g_i(z^t)$ also has a slope equal to
$t$. By Lemma~\ref{lem:ramification_minimal}, the minimal operator of
$g_i(z)$ has a slope equal to $1$. This contradicts Step~2.

\subsection*{Step 4}

\textit{The minimal operator annihilating a $G$-function is regular
singular at $0$.}

\medskip

By~(b), it is enough to prove that its Newton--Ramis polygon has no
positive slope. Suppose, on the contrary, that the minimal operator
$L$ of a $G$-function $g$ has a positive slope $p/q$, where $p$ and
$q$ are positive integers. The series $g(z^q)$ is again a
$G$-function. By Lemma~\ref{lem:ramification_minimal}, the
Newton--Ramis polygon of its minimal operator has the positive
integral slope $p$. This contradicts Step~3. Thus $\mathcal N(L)$ has
no positive slope, and $L$ is regular singular at $0$.

\subsection*{Step 5}

\textit{The non-zero slopes of the Newton--Ramis polygon of the minimal
operator of an $E$-function belong to $[-1,0)$. In particular, this
operator is regular singular at $0$.}

\medskip

Let $f$ be an $E$-function and let $g$ be its associated $G$-function. Let $L$ denote the minimal
operator of the $G$-function $zg(z)$. Let $\Psi$ be the operator obtained from $L$ by the change
of variable $z\mapsto z^{-1}$. Thus, $\Psi$ annihilates
$z^{-1}g(z^{-1})$. By~(f), the operator
$
\Phi:=\overline{\mathcal F}(\Psi)
$
annihilates $f$. 
Let $t\neq0$ be a finite slope of $\mathcal N(\Phi)$. By~(c),
$
-\frac{t}{t+1}
$
is a slope of $\mathcal N(\Psi)$. After the change of variable
$z\mapsto z^{-1}$, property~(d) shows that
$
\frac{t}{t+1}
$
is a slope of $\mathcal N(L)$. 
By Step~4, the minimal operator $L$ has no positive slope. Hence
$
\frac{t}{t+1}\leq0,
$
unless $t=-1$, in which case the quotient is infinite. It follows
that
$
-1\leq t<0.
$
The minimal operator of $f$ is a right factor of $\Phi$. Property~(a)
therefore shows that all of its finite non-zero slopes also belong to
$[-1,0)$. Moreover, by Step~1, the minimal operator of $f$ is a right
factor of an operator having no infinite slope, and hence it has no
infinite slope either. In particular, it has no positive slope and is
regular singular at $0$.

\medskip

We are now ready to conclude the proof of
Theorem~\ref{thm:Andre_consequence}. The finite non-zero slopes of the
local Newton polygon at infinity are the opposites of the negative
Newton--Ramis slopes. Hence, Step~5 implies that the finite slopes of
the local Newton polygon at infinity of the minimal operator of an
$E$-function belong to $[0,1]$, as claimed. It only remains to check
that $\infty$ cannot be regular singular when $f$ is non-polynomial.
An entire solution of a differential equation that is regular singular
at infinity has at most polynomial growth at infinity. By Cauchy's
estimates, such an entire function is a polynomial. Thus the minimal
differential equation of a non-polynomial $E$-function is irregular at
infinity. This completes the proof of
Theorem~\ref{thm:Andre_consequence}. \qed 

\subsection*{Declaration on the use of AI}
All new mathematical ideas and results presented in this article
are entirely due to the authors. AI tools were used solely
to check spelling and identify typographical errors in
mathematical expressions.


\end{document}